\documentclass[twocolumn,secthm]{autart}

\usepackage{amssymb,amsmath,amsfonts}  % assumes amsmath package installed
\usepackage{graphicx}
\usepackage{epstopdf}

\usepackage{mathtools} 
\usepackage[ruled,norelsize,algo2e]{algorithm2e}
\usepackage{comment}
\usepackage{bm}

\usepackage{tikz}
\usetikzlibrary{positioning}
\usepackage{setspace}
\usepackage{lineno}
\usepackage{cite}
\usepackage{url}
\usepackage{booktabs}
\usepackage{multirow}
\usepackage{dsfont}
\usepackage{subcaption}
\usepackage{enumitem,kantlipsum}

\usepackage{caption}% http://ctan.org/pkg/caption
\allowdisplaybreaks

\newcommand{\defeq}{:=}
\newcommand{\dd}{\text{d}}

\newcommand{\innerproduct}[2]{\langle#1,#2\rangle}
\newcommand{\norm}[1]{\left\lVert#1\right\rVert}

\DeclareMathOperator*{\esssup}{ess\,sup}
\newcommand{\trace}[1]{\text{Tr}(#1)}

\newcommand{\costeval}[2]{J_{#1}(#2)}
\newcommand{\optcosteval}[2]{J^*_{#1}(#2)}

\newcommand{\revision}[1]{{\color{black}#1}}

\begin{document}

\begin{frontmatter}
\date{01 November 2021}
\title{Optimal guidance and estimation of a 2D diffusion-advection process by a team of mobile sensors\thanksref{footnoteinfo}} % Title, preferably not more 

\thanks[footnoteinfo]{This paper was not presented at any IFAC 
meeting. Corresponding author Sheng Cheng Tel. +1 217 979 8496.}

\author[Sheng]{Sheng Cheng}\ead{chengs@illinois.edu},    % Add the 
\author[Derek]{Derek A. Paley}\ead{dpaley@umd.edu}               % e-mail address 
% \author[Baiae]{Publius Maro Vergilius}\ead{vergilius@culture.ir}  % (ead) as shown

\address[Sheng]{University of Illinois Urbana-Champaign}  % Please supply                                              
\address[Derek]{University of Maryland, College Park}             % full addresses
% \address[Baiae]{The White House, Baiae}        % here.

\begin{keyword}                           % Five to ten keywords, 
Infinite-dimensional systems; Multi-agent systems; Modeling for control optimization; Guidance navigation and control. 
% chosen from the IFAC 
\end{keyword}                             % keyword list or with the 
                                          % help of the Automatica 
                                          % keyword wizard

\begin{abstract}                          % Abstract of not more than 200 words.
This paper describes an optimization framework to design guidance for a possibly heterogeneous team of multiple mobile sensors to estimate a spatiotemporal process modeled by a 2D diffusion-advection process. Owing to the abstract linear system representation of the process, we apply the Kalman-Bucy filter for estimation, where the sensors provide linear outputs. We propose an optimization problem that minimizes the sum of the trace of the covariance operator of the Kalman-Bucy filter and a generic mobility cost of the mobile sensors, subject to the sensors' motion modeled by linear dynamics. We establish the existence of a solution to this problem. Moreover, we prove convergence to the exact optimal solution of the approximate optimal solution. That is, when evaluating these two solutions using the original cost function, the difference becomes arbitrarily small as the approximation gets finer. To compute the approximate solution, we use Pontryagin's minimum principle after approximating the infinite-dimensional terms originating from the diffusion-advection process. 
The approximate solution is applied in simulation to analyze how a single mobile sensor's performance depends on two important parameters: sensor noise variance and mobility penalty. We also illustrate the application of the framework to multiple sensors, in particular the performance of a heterogeneous team of sensors.
\end{abstract}

\end{frontmatter}
% \linenumbers

% TODO:
% 1. FIXED: switch the usage of alpha and (a,b) for diffusion coefficient and mobile sensor dynamics, respectively.
% 7. Add a pro-and-con analysis of different criteria of improving estimation performance
% 8. Maybe we need to change the reference of Bucy's paper to  1968
% 9. FIXED: We need to change the a and b used in the sensors' dynamics to alpha and beta
% 10. FIXED: We need to change the arrow errorneously used in equations 4.20 and 4.21.
% 11. NOTE: For the optimal control existence of more general problem setting, see Theorem 23.11 in "Functional Analysis, Calculus of Variations and Optimal Control" by Francis Clarke.
% 12. NOT AN ISSUE: The sentence "factors related to the mobile sensor platforms are integrated in addition to reducing the estimation uncertainty only" may not be clear. Try to rewrite it.
% 13. Fixed Figure 2(a) should remove the unit (m) in the the labels.

% ERROR: search ERROR to locate it.

\section{Introduction}
The modern manufacturing industry has benefited from the advantages of mobile robots for their reliability, economic efficiency, safety, and ease of use. However, the monitoring and control of large-scale spatiotemporal processes, e.g., oil spills, harmful algal blooms, and forest fires, have relied heavily on human operators. These events can pose health threats, cause severe environmental issues, and incur substantial financial costs. In such a situation, mobile robots can replace the human operators and carry out the tasks intelligently and cooperatively so long as they are informed about the environment and guided with proper methods.

Spatiotemporal processes evolve in both space and time and, hence, they can be treated as dynamical systems whose dynamics are modeled by partial differential equations (PDEs).
It is generally impossible to completely identify the state of a PDE-modeled system, also known as a distributed parameter system (DPS), with a finite number of sensors. Hence, an observer for the DPS is necessary. For system-theoretical results on the observability of parabolic PDEs, refer to \cite{Khapalov1994exact,Khapalov1993Continuous}. Early designs of the observer include least-square filtering methods for systems governed by linear \cite{Meditch1971} and nonlinear \cite{Lamont1972state} partial differential equations. 

For a DPS whose infinite-dimensional dynamics is linear with additive Gaussian noise, an infinite-dimensional Kalman-Bucy filter (KF), which is an extension to its finite-dimensional analog \cite{Kalman1960ANA,bucy2005filtering}, can be applied. The infinite-dimensional version of the KF first appeared in \cite{bensoussan1971filtrage}. The key to an infinite-dimensional KF is the evolution of the estimation error's covariance, which is operator-valued and can be solved via Riccati equations whose properties have been discussed in \cite{Burns2015solutions,Burns2015Infinitedimensional}. 
%Specifically, conditions for the existence of Bochner integrable solution of the infinite-dimensional Riccati integral equation have been established.
For numerical approximation and computational issues, approximation results are summarized in \cite{Burns2008Mesh} for the infinite-dimensional algebraic Riccati equations of a linear-quadratic regulator.
However, often the disturbance to a system may not be Gaussian, in which case an $H_2$- or $H_{\infty}$-observer design is favorable \cite{Wang2017spatial,Morris2020Optimal}. 
%For a comparison of various observer designs for the heat equation, see \cite{Afshar2015comparison}.

We categorize the design of an observer for a PDE-modeled system by sensor location, i.e., on the boundary or in the interior of the PDE's domain.
For sensors placed on the boundary, one may design observers based on boundary measurements. The challenge in using boundary measurements is that the output operator is unbounded when characterized in an abstract linear system. A tutorial paper \cite{Emirsjlow2000Boundary} reviews the conditions for the well-posedness of the dual control problem with boundary inputs. One approach to boundary observer design is backstepping \cite{Smyshlyaev2005Backstepping,Karagiannis2019backstepping}, which stabilizes the observer using a Volterra transformation that transforms the original system to a stable target system and passes on the stability to the observer via its inverse transformation. %In \cite{Moura2014adaptive}, an adaptive boundary observer was designed using backstepping to estimate battery state-of-charge and state-of-health. 
Optimization techniques, e.g., a linear-quadratic estimator \cite{Moura2011optimal}, have been proposed for boundary sensors using the method of variation.

For static in-domain sensors, a network can be deployed for estimating a PDE-modeled system. The problem is how to place the sensors to yield effective estimation, which is referred to as sensor placement. The solution is to select the sensor placement that yields optimal performance for certain criteria. The trace of the covariance operator of the KF, which quantifies the uncertainty of the estimate, is a common choice of the cost function to be minimized \cite{Bensoussan1972optimization}. A similar problem is investigated in \cite{Zhang2018sensor}, which establishes the well-posedness of the sensor placement problem and its approximation with the cost function being the trace of the covariance operator. %In \cite{Hintermuller2017optimal}, a sensor placement problem was proposed with a trace-involved term as the cost functional and an integral Riccati equation as the constraint. 
The same criterion has been applied to sensor placement of the Boussinesq equation~\cite{Hu2016sensor}. 
In \cite{privat2015optimal}, a randomized observability constant is minimized by choosing suitable shapes and locations of the sensors.  
% Other criteria, e.g., enhanced observability, optimal state estimation, and robust input-output mapping, are discussed for a parabolic PDE in \cite{Demetriou2004optimization}. 
Note that observability may not be a useful criterion \cite{Zhang2018sensor} based on the analysis of its dual, controllability. Specifically, maximizing the minimum eigenvalue of the controllability gramian is not useful (for actuator placement) \cite{morris2015comparison}, because the minimum eigenvalue approaches zero as the dimension of approximation increases in computation.
Geometric approaches can also be applied to sensor placement, e.g., using the centroidal Voronoi tessellation \cite{Demetriou2018using} and combining the transfer-function model with geometric rules \cite{Veldman2020sensor}.
% Geometric approaches can also be applied to sensor placement. \cite{Demetriou2018using} proposes a scheme that places sensors using the centroidal Voronoi tessellation of the kernel of the observer gain of a parabolic PDE. \cite{Veldman2020sensor} proposes a method that combines the transfer-function model and geometric rules to design sensor and actuator locations for which high-gain and low-gain proportional feedback control can reduce the influence of pointwise disturbances.

A variation of the sensor placement problem includes mobility of the sensors. In this scenario, a guidance policy is necessary to take advantage of the additional degree of freedom induced by mobility, which also makes the problem more complicated by introducing the dynamics of the mobile sensors. One may design sensor guidance using Lyapunov-based methods, where the guidance is constructed to make the derivative of the Lyapunov function negative. The Lyapunov function is usually designed to contain (quadratic) terms involving the PDE state and the sensor guidance \cite{Demetriou2010guidance,Demetriou2009estimation}. The Lyapunov-based guidance can further be used in monitoring a hazardous environment where the regions of high information density reduce sensor longevity. Such guidance is combined with a switching policy to balance the conflicting needs of information collection and sensor life span \cite{Demetriou2019incorporating}. %A similar approach uses the gradient of estimation error to guide sensors to the region that has large estimation error \cite{Demetriou2016domain}.

Optimization can also be applied to design sensor guidance, where the guidance (or the trajectory) of the sensor is selected to minimize a cost function. An early work \cite{Carotenuto1987Optimal} proposes an optimization problem that minimizes the weighted sum of the guidance effort for steering a sensor and the mean-square estimation error at a terminal time. In \cite{Demetriou2016gain}, the sensors are guided to the location that yields the maximum value of the estimation kernel. In \cite{imsland2017partially}, receding horizon guidance is proposed to find the sensor path that maximizes mutual information between sensor measurements and the predicted state of the PDE.

This paper proposes an optimization framework that designs guidance for a team of mobile sensors to estimate a 2D diffusion-advection process using a centralized KF. The cost to be minimized is the sum of two terms: the trace of the covariance operator of the KF, called the \textit{uncertainty cost}, which quantifies the uncertainty of the estimation error and the \textit{mobility cost} associated with the sensors' motion. The covariance operator of the KF, which is the solution of an operator-valued Riccati equation, has been studied in \cite{Burns2015solutions}. Specifically, conditions for the existence of Bochner integrable solutions (with values in the Schatten $p$-class) of an operator-valued Riccati equation are established. The Bochner integrable solutions yield simple numerical quadratures for computation of the covariance operator, which is demonstrated in sensor placement \cite{Burns2015solutions} and sensor trajectory planning \cite{Burns2015Infinitedimensional}. Both problems minimize the trace of a weighted covariance operator, whereas the latter has the sensors' dynamics as the constraint.

\revision{In our formulation, factors related to the mobile sensor platforms are integrated in addition to reducing the estimation uncertainty only \cite{Burns2015Infinitedimensional}. Specifically, the integration is reflected by the mobility cost, which 
%includes generic cost functions of the sensors' mobility which 
can be interpreted as the penalty associated with motion. Our formulation supports various types of cost functions for evaluating the cost or penalty induced from mobility. In prior work, the mobility cost is either limited to a quadratic guidance effort \cite{Carotenuto1987Optimal} or cast as a general function of the guidance without detailed discussion \cite{burns2010optimal,Bensoussan1972optimization}. Furthermore, the formulation permits treating the proposed problem as an optimal control problem, where we show both the uncertainty cost and the mobility cost are continuous mappings of the sensors' guidance. We use the techniques for the existence of an optimal control \cite[Theorem~6.1.4]{werner2013optimization} to establish the existence of a solution to our problem. %(Theorem~\ref{thm: existence of a solution of minVarGuidanceProblem})
}

To compute an optimal solution, approximations of the infinite-dimensional terms are necessary. Our treatment of the proposed problem (and its approximation) permits the application of a two-point boundary value problem derived using Pontryagin's minimum principle. After restricting the admissible guidance functions to a stringent set (with a reasonable physical interpretation), we establish convergence to an exact optimal solution of an approximate optimal solution, %(Theorem~\ref{thm: convergence of approximate solution minVarEst}) in the sense that their costs evaluated by the original problem get arbitrarily close as the dimension of approximation increases.
i.e., the cost difference between the original and approximate solutions becomes arbitrarily small as the dimension of approximation increases.
The convergence result justifies the use of the approximation and affirms that the performance of an approximate solution is arbitrarily close to the performance of an exact solution.  

We implement the solution method numerically in simulations to evaluate and analyze the performance of a single sensor and multiple sensors. The flow field that yields advection is set to drift the sensor platforms under realistic conditions. It has been observed that the flow field is leveraged by the optimal guidance to reduce the mobility cost in both cases of a single sensor and a team of homogeneous sensors. Simulations with a team of heterogeneous sensors suggest that such a team can reduce the cost of investment with only a minor degradation in the overall performance.

The contributions of this paper are threefold. First, we propose a guidance policy for a team of mobile sensors to estimate a 2D diffusion-advection process by solving a newly formulated optimization problem that minimizes the trace of the covariance operator plus a generic cost of the sensors' motion subject to sensor platform dynamics. The resulting guidance minimizes the estimation uncertainty while taking care of the sensor-oriented concerns. The formulation with a generic mobility cost applies to a wide range of applications, e.g., accumulated exposure to hazards, guidance effort, or distance to terminal rendezvous locations. The mobility cost is incorporated into the design process, unlike prior work that minimized the uncertainty cost only \cite{Burns2015Infinitedimensional}. Second, we establish realistic conditions for the existence of a solution to the proposed problem. Third, we establish conditions for convergence to an exact optimal solution of the approximate optimal solution. When evaluating these two solutions using the original cost function, the difference gets arbitrarily small as the approximation gets finer. 
%; and (3) solves the approximated problem using the optimal control method.}
%analyzes via simulation the impact on the performance of the proposed guidance of sensor noise and mobility penalty. 

The problem studied in this paper is an extension of the optimal guidance problem for a 1D diffusion process \cite{Cheng2020optimalguidance}. It is the dual problem of the one studied in \cite{Cheng2020optimalInReivew}, which simultaneously designs guidance and actuation of a team of mobile actuators to control a 2D diffusion-advection process.

The remainder of the paper is organized as follows. Section~II introduces the dynamics of the sensors, an abstract linear-system representation of the diffusion-advection process, the measurement model, and the infinite-dimensional Kalman-Bucy filter. Section~III states the problem formulation and establishes conditions for the existence of a solution to the problem. % and proves the existence of a solution.
Section~IV proves the convergence to the exact optimal solution of the approximate optimal solution and introduces a solution method to obtain optimal guidance. Section~V includes the simulation results of multiple parameter studies: a single sensor, a team of homogeneous sensors, and a team of heterogeneous sensors. Section~VI summarizes the paper and discusses ongoing work.

%\singlespacing
% \iffalse This is not. \fi This is typeset again.

\section{Background}
\subsection{Notation and terminology}
The paper adopts the following notation. The symbols $\mathbb{R}$ and $\mathbb{R}^+$ denote the set of real numbers and the set of nonnegative real numbers, respectively. The $n$-nary Cartesian power of a set $M$ is denoted by $M^n$. The notation $X_1 \hookrightarrow X_2$ means that the space $X_1$ is densely and continuously embedded in $X_2$. The norm in a finite- and infinite-dimensional space is denoted by $|\cdot|$ and $\norm{\cdot}$, respectively, with subscripts indicating its type.
The space of all bounded linear operators from space $X$ to space $Y$ is denoted by $\mathcal{L}(X;Y)$ or $\mathcal{L}(X)$ if $Y=X$. We define the space of continuous mappings by $C(I;X) = \{F: I \rightarrow X \text{ such that } t\mapsto F(t) \text{ is continuous in } \norm{\cdot}_X \}$ with the sup norm $\norm{F(\cdot)}_{C(I;X)} = \sup_{t \in I}\norm{F(t)}_X$. 
For a Hilbert space $\mathcal{H}$ equipped with inner product $\innerproduct{\cdot}{\cdot}$ and $\phi_1, \phi_2 \in \mathcal{H}$, define $\phi_1 \circ \phi_2 \in \mathcal{L}(\mathcal{H})$ by $(\phi_1 \circ \phi_2) \psi = \phi_1 \innerproduct{\phi_2}{\psi}$ for all $\psi \in \mathcal{H}$%\cite[Definition~2.1]{Zhang2018sensor}.
. The superscript~${}^*$ denotes an optimal variable, whereas ${}^{\star}$ denotes the adjoint of a linear operator. The transpose of a matrix $A$ is denoted by $A^\top$. An $n \times n$-dimensional diagonal matrix with elements of vector $[a_1,a_2,\dots,a_n]$ on the main diagonal is denoted by $\text{diag}(a_1,a_2,\dots,a_n)$. %The derivative of a function $f$ evaluated at $x$ is denoted by $\dot{f}(x)$.The trace of an operator $\Pi$ and a square matrix $P$ is denoted by $\trace{\Pi}$ and $\traceMatrix{P}$, respectively. 
The $i$th element of a vector $v$ is $[v]_i$. We follow the terminology of \cite{Cheng2020optimalInReivew}: guidance refers to steering the dynamics of the mobile sensor. For an optimization problem
\begin{equation}
\begin{aligned}
    & \underset{x}{\text{minimize}} && J(x) \\ %\nonumber \\
    & \text{subject to} & & x \in C,
\end{aligned}
\tag{P0}
\label{prob: notation example problem}
\end{equation}
we use $\costeval{\eqref{prob: notation example problem}}{x}$ to denote the cost function of \eqref{prob: notation example problem} evaluated at $x$. Specifically, $\optcosteval{\eqref{prob: notation example problem}}{x^*}$ indicates that the optimal value of \eqref{prob: notation example problem} is attained when the cost function is evaluated at an optimal solution $x^*$.

\subsection{Dynamics of the mobile sensors} \label{sebsec: mobile sensor dynamcis}
Assume each of the $m_s$ mobile sensors has the linear dynamics
\begin{equation}\label{eq: agentwise linear dynamics}
    \dot{\zeta}_i(t) = \alpha_i \zeta_i(t) + \beta_i p_i(t), \quad \zeta_i(0) = \zeta_{i,0},
\end{equation}
where $\zeta_i(t) \in \mathbb{R}^{n}$ and $p_i(t) \in P_i \subset \mathbb{R}^{m}$ are the state and the guidance of sensor $i$ at time $t$, respectively. The state $\zeta_i$ contains the 2D location of sensor $i$ and hence $n \geq 2$. Assume that system \eqref{eq: agentwise linear dynamics} is controllable. %The first two elements of $\zeta_i(t)$ are the horizontal and vertical position, $x_i(t)$ and $y_i(t)$, of the sensor in the 2D domain. 
One special case of \eqref{eq: agentwise linear dynamics} would be a single integrator, where $\zeta_i(t) \in \mathbb{R}^2$ is the location, $p_i(t) \in \mathbb{R}^2$ is the velocity command, and $\alpha_i$ and $\beta_i$ are the zero matrix and identity matrix, respectively.

% The last two elements are the horizontal and vertical velocities $\dot{x}_i(t)$ and $\dot{y}_i(t)$. The guidance $p_i(t)$ contains $p_{i,x}(t) \in \mathbb{R}$ along the horizontal axis and $p_{i,y}(t) \in \mathbb{R}$ along the vertical axis.

For conciseness, we concatenate the states and guidance of all $m_s$ sensors, respectively, and use one dynamical equation to describe the dynamics of all agents:
\begin{equation}\label{eq: general dynamics of the mobile sensor}
    \dot{\zeta}(t) = \alpha \zeta(t) + \beta p(t), \quad \zeta(0) = \zeta_0,
\end{equation}
where matrices $\alpha$ and $\beta$ are assembled from $\alpha_i$ and $\beta_i$
for $i \in \{1,2,\dots,m_s\}$, respectively, and are consistent with the concatenation for $\zeta$ and $p$.
% where
% \begin{align}
%     \zeta^\top(t) = & \begin{bmatrix}
%     \zeta_1^\top(t) & \zeta_2^\top(t) & \dots & \zeta_{m_s}^\top(t)
%     \end{bmatrix}, \nonumber \\
%     \zeta_0^\top = & \begin{bmatrix}
%     \zeta_{1,0}^\top & \zeta_{2,0}^\top & \dots & \zeta_{m_s,0}^\top
%     \end{bmatrix}, \nonumber \\
%     p^\top(t) = & \begin{bmatrix}
%         p_1^\top(t) & p_2^\top(t) & \dots & p_{m_s}^\top(t)
%     \end{bmatrix}, \nonumber \\
%     \alpha = & \text{diag}(\alpha_1,\alpha_2,\dots,\alpha_{m_s}), \nonumber \\
%     \beta = & \text{diag}(\beta_1,\beta_2\dots, \beta_{m_s}). \nonumber
% \end{align}
The controllability of the concatenated system \eqref{eq: general dynamics of the mobile sensor} inherits that of each individual system \eqref{eq: agentwise linear dynamics}. With a slight abuse of notation, we use $n$ for the dimension of $\zeta(t)$ and $m$ for the dimension of $p(t)$. Define the admissible set of guidance $P \defeq P_1 \times P_2 \times \dots \times P_{m_s}$ such that $p(t) \in P$ for $t \in [0,t_f]$. %Let $M \in \mathbb{R}^{2 m_s \times n}$ be a matrix such that $M \zeta(t)$ is a vector of locations of the sensors.

\subsection{Abstract linear system}

Consider the following inhomogeneous diffusion-advection equation over a smooth and bounded 2D spatial domain $\Omega$:
\begin{align}
    \frac{\partial z(x,y,t)}{\partial t} = & \ a \nabla^2 z(x,y,t) - \mathbf{v} \cdot \nabla z(x,y,t)  \nonumber \\
    & + w(x,y,t) \label{eq: dynamics of the diffusion advection}\\
    z(x,y,0) =& \ \hat{z}(x,y,0) + w_0(x,y) \\
    z(x,y,t)|_{(x,y) \in \partial \Omega} = & \ 0,
\end{align}
where $(x,y) \in \Omega$, $t \in [0,t_f]$, $a>0$ is the diffusion coefficient, and $\mathbf{v} \in \mathbb{R}^2$ is the flow that yields advection. The initial condition $z(\cdot,\cdot,0)$ is perturbed around its nominal value $\hat{z}(\cdot,\cdot,0)$ by an additive zero-mean Gaussian noise $w_0(\cdot,\cdot)$. The dynamics \eqref{eq: dynamics of the diffusion advection} is subject to an additive zero-mean Gaussian noise $w(\cdot,\cdot,t)$ with variance $Q(t)$, which is nonnegative and self-adjoint. The state noise $w(\cdot,\cdot,t)$ and initial noise $w_0(\cdot,\cdot)$ are mutually independent for all~$t$. % In Ruth Curtain's paper 'A Survey of Infinite Dimensional Filtering', p. 406, for a Wiener process $w(t)$ with incremental covariance $W$ has the following property:
% \mathbb{E}(||w(t)-w(s)||^2) = Tr(W)|t-s|
% where W is positive, self-adjoint, and nuclear (trace class).

% {\color{blue} The dynamics of the mobile sensors follow from the Automatic settings.}
% Assume linear dynamics for sensor $i$, whose location $\zeta_i$ can be controlled via guidance $p_i$ such that
% \begin{align}\label{eq: vehicle dynamics}
% \dot{\zeta}_i(t) = a_i \zeta_i(t) + b_i p_i(t), \ \zeta_i(0) = \zeta_{i0}.
% \end{align}
% We use $a = \text{diag}(a_1,a_2,\dots,a_m)$, $b = \text{diag}(b_1,b_2,\dots,b_m)$, $p = (p_1,p_2,\dots,p_m)^\top$, and $\zeta = (\zeta_1,\zeta_2,\dots,\zeta_m)^\top$ for conciseness. 

% The set of guidance $U$ is defined as $U=\{ p:  p$ is measurable, %pointwise bounded $|p(t)| \leq p_{\text{max}}$, $\forall t \in I$, 
% uniformly bounded by $p_{\text{max}}>0$, and Lipschitz continuous $|p(t_1) - p(t_2)| \leq c_0 |t_1-t_2|$ for $t_1,t_2 \in I \}$.
% By the Arzel\`{a}–Ascoli theorem, the set $U$ is compact if the distance defined on $U$ is the max-norm such that $d(f,g) = \underset{t \in I}{\max}|f(t)-g(t)|$, for $f,g \in U$.
% This definition is one of several ways to construct a compact set of guidance. We choose this definition because of its clear physical interpretation: the bound $p_{\text{max}}$ represents the maximum speed of a vehicle, whereas the Lipschitz coefficient $c_0$ represents the maximum acceleration of the vehicle.

For simplicity, represent the PDE state in \eqref{eq: dynamics of the diffusion advection} by an abstract linear system whose state variable $\mathcal{Z}(t)$ represents $z(\cdot,\cdot,t)$ at time $t$, such that 
\begin{equation}\label{eq: abstract linear system}
\dot{\mathcal{Z}}(t) = \mathcal{A} \mathcal{Z}(t)  + w(t), \quad \mathcal{Z}(0) = \hat{\mathcal{Z}}_0 + w_0 ,
\end{equation}
where $\mathcal{Z}$ belongs to a Hilbert space $\mathcal{H}$ with inner product $\innerproduct{\cdot}{\cdot}$ and induced norm $\norm{\cdot}_{\mathcal{H}}$. Here, the variable $\mathcal{Z}$ is the state of the DPS and space $\mathcal{H}=L^2(\Omega)$ is the state space. 
The operator $\mathcal{A}$ is defined as $\mathcal{A} \psi = a \nabla^2 \psi - \mathbf{v}\cdot \nabla \psi$ with $\psi \in \text{Dom}(\mathcal{A}) = \{ {\psi \in H_0^1(\Omega),} \ {\nabla^2 \psi \in L^2(\Omega)} \} = H^2(\Omega) \cap H_0^1(\Omega)$ \cite{demetriou2012adaptive}. % The operator $\mathcal{D}(\cdot) \in L^2(I,\mathcal{L}(\mathbb{R},\mathcal{H}))$ is the operator version of $D(\cdot,\cdot)$ in~\eqref{eq: dynamics of the diffusion advection}. 
%TO BE COMMENTED: The belonging can be verified by
%\begin{equation*}
%	\int_I \sqrt{\int_{\Omega} (\mathcal{D}(t)a)^2(x) \dd x } \dd t = a \int_I \norm{D(\cdot,t)}_{L^2(\Omega)} \dd t < \infty.
%\end{equation*}
Let $\mathcal{S}(\cdot)$ be the strongly continuous semigroup generated by $\mathcal{A}$.

The measurement by sensor $i$ depends on the sensor's location such that
\begin{equation}
    y_i(t) = \mathcal{C}_i^{\star}(M_0\zeta_i(t),t) \mathcal{Z}(t) +  v_i(t),
\end{equation}
where $M_0$ is a matrix with appropriate dimension such that $M_0 \zeta_i(t) \in \mathbb{R}^2$ is the location of sensor $i$ and $\mathcal{C}_i^{\star} (M_0 \zeta_i(t),t) \in \mathcal{L}(\mathcal{H};\mathbb{R})$ is the output operator that has an integral kernel $\mathcal{C}_i(M_0 \zeta_i(t),t) \in L^2(\Omega)$ such that
\begin{equation*}
    \mathcal{C}_i^{\star}(M_0 \zeta_i(t),t) \phi = \iint_{\Omega} \mathcal{C}_i(M_0 \zeta_i(t),t)(x,y) \phi(x,y) \dd x \dd y
\end{equation*}
for $\phi \in \mathcal{H}$. Additive zero-mean Gaussian noise $v_i(t)$ with variance $\sigma_i^2$ is included in the measurement.

The measurement can have many types, e.g., pointwise \cite{Carotenuto1987Optimal,Khapalov1993Continuous,Khapalov1994exact}, interval integral \cite{Burns2015Infinitedimensional,Demetriou2009estimation}, interval average \cite{Demetriou2009scheduling}, and Gaussian-type kernel \cite{Burns2015Infinitedimensional}. Later in the simulation section, we will use a time-invariant kernel given by the square-shaped average 
\begin{multline}\label{eq: examplery measurement kernel function}
    \mathcal{C}_i(M_0 \zeta_i(t))(x,y) \\ =  \left\{ 
    \begin{aligned}
        \frac{1}{4r_i^2} & \ \text{ if } \begin{bsmallmatrix}
        x \\y
        \end{bsmallmatrix} - M_0 \zeta_i(t) \in [-r_i,r_i] \times [-r_i,r_i] 
        \\
        0 & \ \text{ otherwise,}
    \end{aligned}\right.
\end{multline}
where $2r_i$ is the length of the side of the square at time $t$.

The measurements $y(t) \in \mathbb{R}^{m_s}$ of all sensors are compactly written as
\begin{equation}\label{eq: abstract measurement function}
    y(t)  = \mathcal{C}^{\star}(M\zeta(t),t) \mathcal{Z}(t) +  v(t),
\end{equation}
where $\mathcal{C}^{\star}(M\zeta(t),t)$ is an operator-valued vector $\mathcal{C}^{\star}(M\zeta(t),t)  \defeq [\mathcal{C}_1^{\star}(M_0\zeta_1(t),t),\dots,\mathcal{C}_{m_s}^{\star}(M_0\zeta_{m_s}(t),t) ]^\top$
% \begin{multline}
%     \mathcal{C}^{\star}(M\zeta(t),t)  \defeq \\ \left[\mathcal{C}_1^{\star}(M_0\zeta_1(t),t),\mathcal{C}_2^{\star}(M_0\zeta_2(t),t),\dots,\mathcal{C}_{m_s}^{\star}(M_0\zeta_{m_s}(t),t) \right]^\top \nonumber
% \end{multline}
and $M \in \mathbb{R}^{2 m_s \times n}$ is a matrix such that $M \zeta(t)$ is a vector of locations of the sensors, i.e., $(M \zeta(t))^\top = [(M_0\zeta_1(t))^\top, \dots , (M_0\zeta_{m_s}(t))^\top]$.
% \begin{equation*}
%     (M \zeta(t))^\top = \begin{bmatrix}
%         (M_0\zeta_1(t))^\top & \dots & (M_0\zeta_{m_s}(t))^\top
%     \end{bmatrix}.
% \end{equation*}
We sometimes use $\mathcal{C}(t)$ for brevity instead of $\mathcal{C}(\zeta(t),t)$, because the sensor state $\zeta(t)$ is a function of $t$.
The measurement noise $v(t)$ is a zero-mean Gaussian vector with covariance $R \defeq \text{diag}(\sigma_1^2,\sigma_2^2,\dots,\sigma_{m_s}^2)$. Assume the noise $w_0(\cdot,\cdot), \ w(\cdot,\cdot,t)$, and $v(t)$ are mutually independent for all $t$. 

\subsection{Kalman-Bucy filter}

Analogous to a finite-dimensional linear system, the infinite-dimensional linear system \eqref{eq: abstract linear system} and \eqref{eq: abstract measurement function} admits a Kalman-Bucy filter (KF). For the derivation of the KF of an abstract linear system, one may refer to \cite{Meditch1971,omatu1989distributed}. %[Chapter 7.3 and 7.4] of omatu1989distributed, to be more specific
The estimation $\hat{\mathcal{Z}}(t)$ of the state $\mathcal{Z}(t)$ can be updated from the measurement $y(t)$ by
\begin{align}
\dot{\hat{\mathcal{Z}}}(t) & = \mathcal{A}\hat{\mathcal{Z}}(t) + \Pi(t) \mathcal{C}(t) R^{-1}(y(t)-\hat{y}(t)),\label{eq: dynamics of estimated state} \\
\hat{y}(t) & =  \mathcal{C}^{\star}(t) \hat{\mathcal{Z}}(t),
\end{align}
with initial condition $\hat{\mathcal{Z}}(0) \defeq \hat{\mathcal{Z}}_0$. The predicted observation of the estimated system is denoted by $\hat{y}(t)$. The covariance operator of the estimation error $\Pi(t) \defeq \mathbb{E}[(\mathcal{Z}(t) - \hat{\mathcal{Z}}(t)) \circ (\mathcal{Z}(t) - \hat{\mathcal{Z}}(t))]$ satisfies the following operator Riccati equation:
\begin{multline}
	\dot{\Pi}(t)  =  \mathcal{A} \Pi(t) + \Pi(t) \mathcal{A}^{\star} +Q(t) \\ - \Pi(t) \bar{\mathcal{C}}\bar{\mathcal{C}}^{\star}(t) \Pi(t), \label{eq: Riccati equation in the weak form}
\end{multline}
where $\bar{\mathcal{C}} \bar{\mathcal{C}}^{\star}(t)$ is a compact representation of $\mathcal{C}(t) R^{-1} \mathcal{C}^{\star}(t)$. The initial condition $\Pi(0)$ is given as the covariance operator $\Pi_0$ of the initial estimation error $\Pi_0 \defeq \mathbb{E}[(\mathcal{Z}(0) - \hat{\mathcal{Z}}(0)) \circ (\mathcal{Z}(0) - \hat{\mathcal{Z}}(0))]$ \cite{Zhang2018sensor}, which is the variance of $w_0$ and is nonnegative and self-adjoint. % see equation (7) in \cite{Zhang2018sensor}

% \begin{defn}[\cite{Burns2015solutions}]
\begin{defn}[\protect{\cite[Definition~1.2]{Burns2015solutions}}]
The trace operator $\trace{\cdot}: \mathcal{L}(\mathcal{H}) \rightarrow \mathbb{R}$ is defined as $\trace{\mathcal{G}} = \sum_{i=1}^{\infty}\innerproduct{\phi_i}{\mathcal{G} \phi_i} $ for nonnegative $\mathcal{G} \in \mathcal{L}(\mathcal{H})$, where $\{\phi_i \}_{i=1}^{\infty}$ is an arbitrary orthonormal basis that spans $\mathcal{H}$.
\end{defn}
Note that $\trace{\mathcal{G}}$ is independent of the choice of the orthonormal basis.

\begin{defn}[\protect{\cite[Definition~3.2]{Burns2015Infinitedimensional}}]
	Let $\mathbb{H}$ be a separable complex Hilbert space. For $1 \leq q < \infty$, let $\mathcal{J}_q(\mathbb{H})$ denote the set of all bounded operators $\mathcal{L}(\mathbb{H})$ such that $\trace{[A]^q} < \infty$, where $[A] \defeq \sqrt{A^{\star}A}$. If $A \in \mathcal{J}_q(\mathbb{H})$, then the $\mathcal{J}_q$-norm of $A$ is defined as $\norm{A}_{\mathcal{J}_q(\mathbb{H})} \defeq (\trace{[A]^q})^{1/q}$.
\end{defn}

The classes $\mathcal{J}_1(\mathcal{H})$ and $\mathcal{J}_2(\mathcal{H})$ are known as the space of trace-class operators and the space of Hilbert-Schmidt operators, respectively. Note that a continuous embedding ${\mathcal{J}_{q_1}(\mathcal{H}) \hookrightarrow \mathcal{J}_{q_2}(\mathcal{H})}$ holds if $1 \leq q_1 < q_2 \leq \infty$ \cite{Burns2015solutions}. In other words, if $A \in \mathcal{J}_{q_1}(\mathcal{H})$, then $A \in \mathcal{J}_{q_2}(\mathcal{H})$ and ${\norm{A}_{\mathcal{J}_{q_2}(\mathcal{H})} \leq \norm{A}_{\mathcal{J}_{q_1}(\mathcal{H})}}$.%\cite{Burns2015Infinitedimensional}.

Consider the following assumptions with $1 \leq q < \infty$:
\begin{enumerate}
    \item[(A1)] $\Pi_0 \in \mathcal{J}_q(\mathcal{H})$ and $\Pi_0$ is nonnegative.
    \item[(A2)] $Q(\cdot) \in L^1([0,t_f];\mathcal{J}_q(\mathcal{H}))$ and $ Q(t) $ is nonnegative for all $t \in [0,t_f]$.
    \item[(A3)] $\bar{\mathcal{C}} \bar{\mathcal{C}}^{\star} (\cdot) \in L^{\infty}([0,t_f];\mathcal{L}(\mathcal{H}))$ and $\bar{\mathcal{C}} \bar{\mathcal{C}}^{\star}(t)$ is nonnegative for $t \in [0,t_f]$.
\end{enumerate}

The existence of a mild solution of \eqref{eq: Riccati equation in the weak form} is established in Lemma~\ref{lemma: existence of Riccati mild solution}. The proof is omitted because the lemma follows directly from \cite[Theorem~3.6]{Burns2015solutions}.

\begin{lem}{\cite[Theorem~3.6]{Burns2015solutions}}\label{lemma: existence of Riccati mild solution}
    Let $\mathcal{H}$ be a separable Hilbert space. Suppose assumptions (A1)--(A3) hold. Then, the equation 
    \begin{multline}
        \Pi(t) = \mathcal{S}(t) \Pi_0 \mathcal{S}^{\star}(t) + \int_0^{t} \mathcal{S}(t-\tau) \\ \left( Q(\tau)  - \Pi(\tau) \bar{\mathcal{C}} \bar{\mathcal{C}}^{\star}(\tau) \Pi(\tau) \right) \mathcal{S}^\star(t-\tau) \dd \tau \label{eq: new mild solution of operator Riccati equation}
    \end{multline}
provides a unique mild solution to \eqref{eq: Riccati equation in the weak form} in the space $L^2([0,t_f];\mathcal{J}_{2q}(\mathcal{H}))$. The solution is in $C([0,t_f];\mathcal{J}_q(\mathcal{H}))$ and is pointwise self-adjoint and nonnegative.
%Furthermore, if $Q(\cdot)\in C([0,t_f];\mathcal{J}_q(\mathcal{H}))$ and $\bar{\mathcal{C}} \bar{\mathcal{C}}^{\star}(\cdot) \in C([0,t_f];\mathcal{L}(\mathcal{H}))$, then $\Pi$ is a weak solution to \eqref{eq: Riccati equation in the weak form}.
\end{lem}

The covariance operator $\Pi(t)$ characterizes the uncertainty of the estimation error. % one paragraph after Lemma 1.4 of Burns2015solutions.
The expected value of the squared norm of the estimation error is the trace of the covariance operator $\Pi(t)$ \cite{Burns2015Infinitedimensional,Zhang2018sensor}: $\trace{\Pi(t)} = \mathbb{E}[\lVert \mathcal{Z}(t) - \hat{\mathcal{Z}}(t) \rVert_{\mathcal{H}}^2 ]$.

The following assumption is vital to the main results in this paper:
\begin{enumerate}
    \item[(A4)] For each sensor $i \in \{1,2,\dots, m_s \}$, the kernel of the output operator $\mathcal{C}_i(x,t)$ is continuous with respect to location $x \in \mathbb{R}^2$ \cite[Definition~4.5]{Burns2015Infinitedimensional}. That is, there exists a continuous function $l: \mathbb{R}^+ \rightarrow \mathbb{R}^+$ such that $l(0) = 0$ and $\norm{\mathcal{C}_i(x_1,t) - \mathcal{C}_i(x_2,t)}_{L^2(\Omega)} \leq l(| x_1  -  x_2|_2)$ for all $t \in [0,t_f]$ and all $x_1,x_2 \in \mathbb{R}^2$.
\end{enumerate}

Assumption (A4) is important in that, roughly speaking, it establishes the continuity of the covariance operator with respect to sensor state (Lemma~\ref{lemma: continuity of trace cost wrt sensor trajectory}), which further permits the existence of a solution to the optimization problem proposed in this paper (Theorem~\ref{thm: existence of a solution of minVarGuidanceProblem}), its finite-dimensional approximation (Theorem~\ref{thm: existence of a solution of the approximated problem}), and the convergence to the exact optimal cost of the approximate optimal cost (Theorem~\ref{thm: convergence of approximate solution minVarEst}).

\begin{rem}\label{remark: continuity wrt location}
	The time invariant kernel in \eqref{eq: examplery measurement kernel function} is continuous with respect to location, where $l(u) = (4r_i c_0 u + u^2)^{1/2} /(4r_i^2)$ for sensor $i$ in assumption (A4) for $c_0 > 0$.% that satisfies $c_0 |x|_1 \leq |x|_2$ for $x \in \mathbb{R}^2$.% \cite[p. 197]{Burns2015Infinitedimensional}.
\end{rem}

The sensors' locations determine where the output is measured and, furthermore, how the covariance operator evolves through \eqref{eq: new mild solution of operator Riccati equation}. We characterize this relation by a composite mapping. Since the output operator $\mathcal{C}^\star(\cdot,t)$ is a mapping of the sensors' locations at time $t$, the composite output operator $\bar{\mathcal{C}} \bar{\mathcal{C}}^\star (\cdot)$ is a mapping of the sensor state in $[0,t_f]$ and so is $\Pi(\cdot)$ by \eqref{eq: new mild solution of operator Riccati equation}, although the sensor state is not explicitly reflected in the notation of the latter two mappings. Hence, we can define the uncertainty cost $ \int_0^{t_f} \trace{\Pi(t)} \dd  t$ as a mapping of the sensor state. Let $K:C([0,t_f];\mathbb{R}^n) \rightarrow \mathbb{R}^+$ such that $K(\zeta) \defeq \int_0^{t_f} \trace{\Pi(t)} \dd  t$. Lemma~\ref{lemma: continuity of trace cost wrt sensor trajectory} below states the continuity of the uncertainty cost with respect to the sensor state. Its proof can be found in the supplementary material.

\begin{lem}\label{lemma: continuity of trace cost wrt sensor trajectory}
    Let assumptions (A1)--(A3) hold with $q = 1$ and $\Pi \in C([0,t_f];\mathcal{J}_1(\mathcal{H}))$ be defined in \eqref{eq: new mild solution of operator Riccati equation}. If assumption (A4) holds, then the mapping $K(\cdot)$ is continuous.
\end{lem}

% \begin{pf*}{Proof}
% See Appendix~\ref{prf: continuity of trace cost wrt sensor trajectory}.
% \qed
% \end{pf*}

% \begin{rem}\label{remark: K() is an implicit mapping}
%     The mapping $K(\zeta) = \int_0^{t_f} \trace{\Pi(t)} \dd  t$ defined in Lemma~\ref{lemma: continuity of trace cost wrt sensor trajectory} is a composite mapping because $\Pi(\cdot)$ is a mapping of the composite output operator $\bar{\mathcal{C}} \bar{\mathcal{C}}^\star(\cdot)$, which is a mapping of the sensor state $\zeta(\cdot)$ by definition.
% \end{rem}

\subsection{Finite-dimensional approximation}\label{subsec: finite dim approximation}
A finite-dimensional approximation of the infinite-dimensional state estimate $\hat{\mathcal{Z}}(t)$ and covariance operator $\Pi(t)$ is necessary for numerical computation.
Consider a finite-dimensional subspace $\mathcal{H}_N \subset \mathcal{H}$ with dimension~$N$. The inner product and norm of $\mathcal{H}_N$ are inherited from that of $\mathcal{H}$. Let $P_N: \mathcal{H} \to \mathcal{H}_N$ denote the orthogonal projection of $\mathcal{H}$ onto $\mathcal{H}_N$. Let $\hat{Z}_N(t) \defeq P_N \hat{\mathcal{Z}}(t)$ and $S_N(t) \defeq P_N \mathcal{S}(t) P_N$ be the finite-dimensional approximations of $\hat{\mathcal{Z}}(t)$ and $\mathcal{S}(t)$, respectively. The approximated version of the estimation \eqref{eq: dynamics of estimated state} is 
\begin{align}
\dot{\hat{Z}}_N(t) & = A_N \hat{Z}_N(t) + \Pi_N(t) C_N(t) R^{-1}(y(t)-\hat{y}_N(t)), \label{eq: dynamics of approximated estimated state} \\
\hat{y}_N(t) & =  C_N^\star(t) \hat{Z}_N(t),
\end{align}
with initial condition $\hat{Z}_N(0)  = P_N \hat{\mathcal{Z}}_0$. The approximations $A_N \in \mathcal{L}(\mathcal{H}_N)$ and $C^{\star}_N(t) \in \mathcal{L}(\mathcal{H}_N;\mathbb{R}^{m_s})$ are of $\mathcal{A}$ and $\mathcal{C}^{\star}(t)$, respectively, and $\Pi_N(t)$ is the finite-dimensional approximation of $\Pi(t)$ such that
\begin{align}\label{eq: mild solution of approximate Riccati covariance}
    \Pi_N(t) = & S_N(t) \Pi_{0,N} S_N^\star(t) + \int_0^t S_N(t - \tau) \big( Q_N(\tau) \nonumber \\ 
    & - \Pi_N(\tau) \bar{C}_N \bar{C}_N^\star (\tau) \Pi_N(\tau) \big) S_N^\star(t - \tau) \dd \tau ,
\end{align}
where $\Pi_{0,N} \defeq P_N \Pi_0 P_N$ and $Q_N(t) \defeq P_N Q(t) P_N$ are approximations of $\Pi_0$ and $Q(t)$, respectively, and $\bar{C}_N \bar{C}_N^\star (\tau)$ is short for $C_N R^{-1} C_N^\star (\tau)$.

If the subspace $\mathcal{H}_N$ is chosen such  that it is spanned by the first $N$ functions of the orthonormal basis $\{\phi_i \}_{i=1}^{\infty}$ that spans $\mathcal{H}$, then $\trace{\Pi_N(t)} = \trace{P_N \Pi(t) P_N} = \sum_{i=1}^N \innerproduct{\phi_i}{\Pi(t) \phi_i}.$
% \begin{equation}
%     \trace{\Pi_N(t)} = \trace{P_N \Pi(t) P_N} = \sum_{i=1}^N \innerproduct{\phi_i}{\Pi(t) \phi_i}.
% \end{equation}
To establish convergence of the approximate covariance operator $\Pi_N(\cdot)$ to the original operator $\Pi(\cdot)$, the following assumptions are made:
\begin{enumerate}
    \item[(A5)] Both $\Pi_0$ and sequence $\{\Pi_{0,N} \}_{N=1}^{\infty}$ are elements of {$\mathcal{J}_q(\mathcal{H})$}. Both $\Pi_0$ and $\Pi_{0,N}$ are nonnegative for all $N \in \mathbb{N}$ and $\norm{\Pi_0-\Pi_{0,N}}_{\mathcal{J}_q(\mathcal{H})} \rightarrow 0$ as $ N \rightarrow \infty$.
    \item[(A6)] Both $Q(\cdot)$ and sequence $\{Q_N(\cdot) \}_{N=1}^{\infty}$ are elements of $L^1([0,t_f];\mathcal{J}_q(\mathcal{H}))$. Both $Q (\tau)$ and $Q_N(\tau)$ are nonnegative for all $\tau \in [0,t_f]$ and all $N \in \mathbb{N}$ and satisfy $\int_0^t \norm{Q(\tau)-Q_N (\tau)}_{\mathcal{J}_q(\mathcal{H})} \dd \tau \rightarrow 0$
        % \begin{equation}
        %     \int_0^t \norm{Q(\tau)-Q_N (\tau)}_{\mathcal{J}_q(\mathcal{H})} \dd \tau \rightarrow 0
        % \end{equation}
    for all $t \in [0,t_f]$ as $N \rightarrow \infty$.
    \item[(A7)] Both $\bar{\mathcal{C}} \bar{\mathcal{C}}^{\star}(\cdot)$ and sequence $\{\bar{C}_N \bar{C}_N^{\star}(\cdot)\}_{N=1}^{\infty}$ are elements of $L^{\infty}([0,t_f];\mathcal{L}(\mathcal{H}))$. Both $\bar{\mathcal{C}} \bar{\mathcal{C}}^{\star}(t) $ and $\bar{C}_N \bar{C}_N^{\star}(t) $ are nonnegative for all $N$ and $t \in [0,t_f]$. And ${\esssup}_{t \in [0,t_f]} \norm{\bar{\mathcal{C}}\bar{\mathcal{C}}^{\star}(t) - \bar{C}_N \bar{C}_N^\star(t)}_{\text{op}} \rightarrow 0$ holds
        % \begin{equation}
        %     \underset{t \in [0,t_f]}{\esssup} \norm{\bar{\mathcal{C}}\bar{\mathcal{C}}^{\star}(t) - \bar{C}_N \bar{C}_N^\star(t)}_{\text{op}} \rightarrow 0
        % \end{equation}
        as $N \rightarrow \infty$ ($\norm{\cdot}_{\text{op}}$ denotes the operator norm).
\end{enumerate}

Note that assumptions (A1), (A2), and (A3) are contained within assumptions (A5), (A6), and (A7), respectively.

The convergence of the approximate covariance operator $\Pi_N(\cdot)$ is stated in the next theorem whose proof is omitted since the theorem is reproduced from \cite[Theorem~3.5]{Burns2015solutions}.

\begin{thm}[\protect{\cite[Theorem~3.5]{Burns2015solutions}}]\label{thm: convergence of riccati covariance operator}
    Suppose $\mathcal{S}(t)$ is a strongly continuous semigroup of linear operators over a Hilbert space $\mathcal{H}$ and that $\{S_N(t)\}$ is a sequence of uniformly continuous semigroup over the same Hilbert space that satisfy
    \begin{equation}\label{eq: primal and dual convergence of the approximation}
        \norm{\mathcal{S}(t) \phi - S_N(t) \phi} \rightarrow 0, \quad \norm{\mathcal{S}^{\star}(t) \phi - S_N^{\star}(t) \phi} \rightarrow 0
    \end{equation}
    % p. 213 of \cite{Burns2015infinte} justifies these two conditions
    as $N \rightarrow \infty$, uniformly in $[0,t_f]$ for each $\phi \in \mathcal{H}$. Suppose assumptions (A5)--(A7) hold. If $\Pi(\cdot) \in C([0,t_f];\mathcal{J}_q(\mathcal{H}))$ is a solution of \eqref{eq: new mild solution of operator Riccati equation} and $\Pi_{N}(\cdot) \in C([0,t_f];\mathcal{J}_q(\mathcal{H}))$ is the sequence of solution of \eqref{eq: mild solution of approximate Riccati covariance}, then $\sup_{t \in [0,t_f]} \norm{\Pi(t)-\Pi_{N}(t)}_{\mathcal{J}_q(\mathcal{H})} \rightarrow 0$ as $N \rightarrow \infty$.
\end{thm}

The following assumption and lemma are related to the continuity with respect to location of the approximate output kernel and the continuity with respect to sensor state of the trace of the approximate covariance operator, which are analogous to assumption (A4) and Lemma~\ref{lemma: continuity of trace cost wrt sensor trajectory}, respectively.

\begin{enumerate}
    \item[(A8)] The approximated input operator $C_{i,N}(x,t)$ is continuous with respect to location $x \in \mathbb{R}^2$, that is, there exists a continuous function $l_N: \mathbb{R}^+ \rightarrow \mathbb{R}^+$ such that $l_N(0) = 0$ and $ \norm{C_{i,N}(x_1,t) - C_{i,N}(x_2,t)}_{L^2(\Omega)} \leq l_N(|x_1-x_2|_2)$ for all $t \in [0,t_f]$, all $x_1,x_2 \in \mathbb{R}^2$, and all $i \in \{1,2,\dots,m_s \}$.
\end{enumerate}

Similar to the mapping $K(\cdot)$ in Lemma~\ref{lemma: continuity of trace cost wrt sensor trajectory}, we can characterize the approximate uncertainty cost as a mapping of the sensor state $\zeta$, where continuity is established in Lemma~\ref{lem: continuity of approximated trace cost wrt sensor trajecotry}, whose proof is in the supplementary material.

\begin{lem}\label{lem: continuity of approximated trace cost wrt sensor trajecotry}
Let assumptions (A5)--(A7) hold and $\Pi_{N}(t)$ be defined as in \eqref{eq: mild solution of approximate Riccati covariance}. If assumption (A8) holds, then the mapping $K_N: C([0,t_f];\mathbb{R}^n) \rightarrow \mathbb{R}^+$ such that $K_N(\zeta) \defeq \int_0^{t_f} \trace{\Pi_N(t)} \dd t$ is continuous.
\end{lem}

% \begin{pf*}{Proof}
% See Appendix~\ref{prf: continuity of approximated trace cost wrt sensor trajecotry}.
% \end{pf*}

% \begin{rem}
%     The mapping $K_N$ defined in Lemma~\ref{lem: continuity of approximated trace cost wrt sensor trajecotry} is a composite mapping of sensor state $\zeta(\cdot)$ for the same reason argued in Remark~\ref{remark: K() is an implicit mapping}.
% \end{rem}

\section{Problem formulation}
We now introduce the formulation of the optimization problem.
Given the dynamics and initial condition \eqref{eq: general dynamics of the mobile sensor} of the sensors, the dynamics of the diffusion-advection process \eqref{eq: abstract linear system}, and the second moment of the initial state noise $w_0(\cdot,\cdot)$, the process noise $w(\cdot,\cdot,t)$, and measurement noise $v(t)$, the problem below yields the optimal guidance for a team of mobile sensors to estimate a 2D diffusion-advection process. 

The cost function consists of two parts: one part accounts for reducing the estimation uncertainty (\textit{uncertainty cost}), and the other accounts for the motion of the sensors (\textit{mobility cost}).
The uncertainty cost is the integral of the trace of the covariance operator $\Pi(\cdot)$ over the horizon $[0,t_f]$, i.e., $\int_0^{t_f} \trace{\Pi(t)} \dd t$.
The mobility cost $J_{\text{m}}(\zeta,p)$ is defined as
\begin{equation}
    J_{\text{m}}(\zeta,p) = \int_0^{t_f} h(\zeta(t),t) + g(p(t),t) \dd  t + h_f(\zeta(t_f)).
\end{equation}
Here, $h : \mathbb{R}^n \times [0,t_f] \rightarrow \mathbb{R}^+$ is a continuous function that characterizes the cost associated with the state of the mobile sensors. For example, a hazardous field can be modeled by $h$, where $h(\zeta(t),t)$ evaluates the exposure of the mobile sensors, which can shorten the sensor's life span. The cost of the guidance is characterized by $g: \mathbb{R}^m \times [0,t_f] \rightarrow \mathbb{R}^+$. For example, a quadratic guidance effort is $g(p(t),t) = p^\top (t) \gamma p(t)$, where $\gamma \in \mathbb{R}^{m \times m}$ is symmetric and positive definite. The guidance cost can address limited onboard resources, like fuel or batteries, by treating $\gamma$ as the penalty coefficient. The terminal state cost $h_f: \mathbb{R}^n \rightarrow \mathbb{R}^+$ evaluates the cost associated with the terminal state. An exemplary scenario is when the sensors are expect to come close to a set of pre-assigned terminal locations $x_f \in \Omega^{m_s}$, where $h_f(\zeta(t_f)) \defeq |M\zeta(t_f)-x_f|_2^2$. 

The motion of the sensors follow from the dynamics \eqref{eq: general dynamics of the mobile sensor}, which constrain the optimization. Denote the admissible set of guidance functions as $\mathcal{P} \defeq  L^2([0,t_f];P)$, where $P$ is the set of admissible guidance (values) defined at the end of Section~\ref{sebsec: mobile sensor dynamcis}.

The optimization problem is formulated as follows:
\begin{equation}\label{prob: min covariance sensor guidance problem}
\begin{aligned}
& \underset{p \in \mathcal{P}}{\text{minimize}} &&  \int_0^{t_f} \trace{\Pi(t)} + h(\zeta(t),t) + g(p(t),t) \dd t \\
& && + h_f(\zeta(t_f)) \\
& \text{subject to}  && \dot{\zeta}(t) = a \zeta(t) + b p(t), \ \zeta(0) = \zeta_0,
\end{aligned}
\tag{P}
\end{equation}
where $\Pi(t)$ is given by \eqref{eq: new mild solution of operator Riccati equation} with a given initial condition $\Pi(0) = \Pi_0$. It suffices to search for guidance $p$ that minimizes the cost of \eqref{prob: min covariance sensor guidance problem}, because the sensor state $\zeta$ is entirely determined by guidance $p$ via the sensor dynamics and the given initial condition $\zeta_0$, which further determines $\Pi(\cdot)$ through \eqref{eq: new mild solution of operator Riccati equation} with a given initial covariance~$\Pi(0)$. 

A special case is considered in \cite{Cheng2020optimalguidance} where only a quadratic guidance effort is considered in the mobility cost. Such a formulation applies to the case of limited onboard resources of each mobile sensor when $\gamma$ is diagonal. It minimizes the Lagrangian function of the optimization problem that minimizes the uncertainty cost subject to the constraints of bounded guidance effort and linear dynamics of the mobile sensors. 

The following three assumptions are necessary for the existence of a solution to problem \eqref{prob: min covariance sensor guidance problem}.
\begin{enumerate}
    %\item The set $X_0 \subset \mathbb{R}^{n_a m_a}$, where $\zeta_0 \in X_0$, is convex and closed.
    \item[(A9)] The set of admissible guidance $P \subset \mathbb{R}^{m}$ is closed and convex.
    \item[(A10)] The mappings $h: \mathbb{R}^{n} \times [0,t_f] \rightarrow \mathbb{R}^+$, $g: \mathbb{R}^{m} \times [0,t_f] \rightarrow \mathbb{R}^+$, and $h_f: \mathbb{R}^{n} \rightarrow \mathbb{R}^+$ are continuous. For every $t \in [0,t_f]$, the function $g(p,t)$ is convex about $p$.
    % The old A9 is removed because it is not very useful in the proof. Instead, we require all the cost functions in the mobility cost are nonnegative.
    \item[(A11)] There exists a constant $d_1 >0$ with $g(p,t) \geq d_1 |p|_2^2$  for all $(p,t) \in P \times [0,t_f]$.
    % \item[(A10)] The operator $\mathcal{B}_i(x,t)$ is continuous with respect to location for $i \in \{1,2,\dots, m_s \}$.
\end{enumerate}

Assumptions (A9)--(A11) are generally met in applications with real vehicles. Assumption (A9) is satisfied when the values of admissible guidance vary along a continuum. The continuity requirement in assumption (A10) on the cost functions $h$, $g$, and $h_f$ is typically satisfied. And the convexity requirement in assumption (A10) and quadratic boundedness from below in assumption (A11) can be met if $g$ is quadratic in $p$, e.g., $g(p,t) = p^\top(t) \gamma(t) p(t)$ for a symmetric and positive definite matrix $\gamma(t)$ (which is continuous with respect to $t$) such that $d_1$ can be chosen as the minimum eigenvalue of $\gamma(t)$ for $t \in [0,t_f]$ (see \cite[Corollary~VI.1.6]{bhatia1996matrix}). 

Theorem~\ref{thm: existence of a solution of minVarGuidanceProblem} below states the existence of a solution to problem \eqref{prob: min covariance sensor guidance problem}, whose proof can be found in Appendix~\ref{prf: existence of a solution of minVarGuidanceProblem}.

 \begin{thm}\label{thm: existence of a solution of minVarGuidanceProblem}
    Consider problem \eqref{prob: min covariance sensor guidance problem} and let assumptions (A1)--(A4) and (A9)--(A11) hold. \revision{Then problem \eqref{prob: min covariance sensor guidance problem} has a solution.}
    % In the first submission
    % Consider problem \eqref{prob: min covariance sensor guidance problem} and let assumptions (A1)--(A4) and (A9)--(A11) hold. If there exists guidance $p_0$ in the admissible set $\mathcal{P}$ such that ${\costeval{\eqref{prob: min covariance sensor guidance problem}}{p_0} < \infty}$, then \eqref{prob: min covariance sensor guidance problem} has a solution.
\end{thm}
% \begin{pf*}{Proof}
%     See Appendix~\ref{prf: existence of a solution of minVarGuidanceProblem}.
%     \qed
% \end{pf*}

We use Pontryagin's minimum principle to characterize an optimal solution of \eqref{prob: min covariance sensor guidance problem}. Consider the Hamiltonian
\begin{multline}\label{eq: Hamiltonian}
H(\zeta(t),p(t),\lambda(t),t) = \trace{\Pi(t)} \\ + h(\zeta(t),t) + g(p(t),t) + \lambda^\top(t) (a \zeta(t) + bp(t) ),
\end{multline}
where $\lambda(t) \in \mathbb{R}^n$ is the costate associated with $\zeta(t)$.
The necessary conditions of (local) optimality are as follows:
\begin{subequations}\label{eq: pontryagin's principle}
\begin{align}
\dot{\zeta}^*(t) = & \ a \zeta^*(t) + b p^* (t),  \label{eq: Pontryagin vehicle dynamics}\\
\zeta^*(0) = & \ \zeta_0, \\
\dot{\lambda}^*(t) = & \ -a^\top \lambda^*(t) - \nabla_{\zeta} h(\zeta^*(t),t) \nonumber \\
& \ - \nabla_{\zeta} \trace{\Pi^*(t)}, \label{eq: Pontryagin costate dynamics}\\	
\lambda^*(t_f) = & \ \nabla_{\zeta} h_f(\zeta^*(t_f)), \\
0 = & \ \nabla_p g(p^*(t),t) + b^\top \lambda^*(t) , \label{eq: Pontryagin optimal control}
\end{align}
\end{subequations}
where $\Pi^*(\cdot)$ is evaluated along the optimal system state $\zeta^*(\cdot)$ and we use the first-order necessary condition $\nabla_{p}H(\zeta^*(t),p^*(t),\lambda^*(t),t) = 0$ in \eqref{eq: Pontryagin optimal control} for $H$ to attain its minimum at $p^*(t)$. The necessary condition \eqref{eq: pontryagin's principle} essentially requires the solution to a two-point boundary value problem, which further requires the derivation of $\nabla_{\zeta} \trace{\Pi^*(t)}$. 
We refer to a similar derivation in \cite{Burns2015Infinitedimensional}, where the gradient of the covariance operator's trace with respect to sensor's guidance is taken.
The $i$th row of $\nabla_{\zeta} \trace{\Pi^*(t)}$, denoted by $[\nabla_{\zeta} \trace{\Pi^*(t)}]_i$, is the partial derivative of $\trace{\Pi^*(t)}$ with respect to the $i$th element of the state $\zeta^*(t)$ for $i \in \{1,2,\dots,n \}$. Since trace is a linear operator, we have 
\begin{equation}\label{eq: partial derivative of trace of Riccati wrt sensor i's location eq 1}
    [\nabla_{\zeta} \trace{\Pi^*(t)}]_i = \frac{\partial \trace{\Pi^*(t)}}{\partial [\zeta(t)]_i} = \trace{\frac{\partial \Pi^*(t)}{ \partial [\zeta(t)]_i}}.
\end{equation}
By the chain rule, \eqref{eq: partial derivative of trace of Riccati wrt sensor i's location eq 1} becomes
\begin{equation}\label{eq: partial derivative of trace of Riccati wrt sensor i's location eq 2}
    \trace{\frac{\partial \Pi^*(t)}{ \partial [\zeta(t)]_i}} = \trace{D_{\bar{\mathcal{C}} \bar{\mathcal{C}}^{\star}(t)} \Pi^*(t) \circ D_{[\zeta(t)]_i} \bar{\mathcal{C}} \bar{\mathcal{C}}^{\star}(t)},
\end{equation}
where $D_{\bar{\mathcal{C}} \bar{\mathcal{C}}^{\star}(t)} \Pi(t) $ is the Fr\'echet derivative of the Riccati operator with respect to the composite output operator $\bar{\mathcal{C}} \bar{\mathcal{C}}^{\star}(t)$ and $D_{[\zeta(t)]_i} \bar{\mathcal{C}} \bar{\mathcal{C}}^{\star}(t)$ is the Fr\'echet derivative of $\bar{\mathcal{C}} \bar{\mathcal{C}}^{\star}(t)$ with respect to $[\zeta(t)]_i$. Denote $D_{\bar{\mathcal{C}} \bar{\mathcal{C}}^{\star}(t)} \Pi^*(t) $ by~$\Lambda(t)$ and, by \cite[Theorem 5.5]{Burns2015Infinitedimensional}, $\Lambda(t)$ is the unique solution to
% integral representation
\begin{align}
	\Lambda h(t) =& \ -\int_0^t \mathcal{S}(t-s) \big((\Lambda h)(s) \bar{\mathcal{C}} \bar{\mathcal{C}}^{\star}(s) \Pi(s) \nonumber \\
	& + \Pi(s)\bar{\mathcal{C}} \bar{\mathcal{C}}^{\star}(s) (\Lambda h)(s) \nonumber \\
	& + \Pi(s) h(s) \Pi(s)\big) \mathcal{S}^{\star}(t-s) \dd s \label{eq: differential equation of the Frechet derivative of the Kalman gain wrt the measurement operator}
\end{align}
%% differential equation representation
%\begin{align}
%\dot{\Lambda} = & -\Lambda(t)\bar{\mathcal{C}}_{\zeta}\bar{\mathcal{C}}_{\zeta}^{\star}(t) \Pi(t) - \Pi(t) \bar{\mathcal{C}}_{\zeta}\bar{\mathcal{C}}_{\zeta}^{\star}(t) \Lambda(t) \nonumber \\
%&- \Pi(t)\Pi(t) + \mathcal{A}\Lambda(t) + \Lambda(t)\mathcal{A}^{\star}, \label{eq: differential equation of the Frechet derivative of the Kalman gain wrt the measurement operator}\\
%\Lambda(0) = & 0.
%\end{align}
with $\Lambda(0) =  0$ for all $h \in C([0,t_f];\mathcal{J}_1(\mathcal{H}))$ and all $t \in [0,t_f]$. The approximated version of problem \eqref{prob: min covariance sensor guidance problem} and the two-point boundary value problem \eqref{eq: pontryagin's principle} will be applied to solve for optimal guidance in Section~\ref{sec: approximation}.

%\subsection{Galerkin approximation}

\section{Solving optimal guidance using approximation}\label{sec: approximation}

% It is generally difficult to derive closed-form solutions for problem \eqref{prob: min covariance sensor guidance problem}. Therefore, numerical approaches will be used to obtain optimal solutions. 
Since the infinite-dimensional terms in \eqref{prob: min covariance sensor guidance problem} have to be approximated for computation as introduced in Section~\ref{subsec: finite dim approximation}, we arrive at the approximated problem:
\begin{equation}\label{prob: approximated min covariance sensor guidance problem}
\begin{aligned}
& \underset{p_N \in \mathcal{P}}{\text{minimize}} &&  \int_0^{t_f} \trace{\Pi_N(t)} + h(\zeta(t),t) + g(p_N(t),t) \dd t\\
& &&  + h_f(\zeta(t_f))\\
& \text{subject to}  && \dot{\zeta}(t) = a \zeta(t) + b p_N(t), \ \zeta(0) = \zeta_0,
\end{aligned}
\tag{AP}
\end{equation}
where $\Pi_N(t)$ is obtained through \eqref{eq: mild solution of approximate Riccati covariance}. It suffices to search for guidance $p_N$, because both the sensor state~$\zeta$ and the approximated estimation covariance $\Pi_N$ are fully determined by the guidance and initial conditions. The existence of a solution of problem \eqref{prob: approximated min covariance sensor guidance problem} is guaranteed in Theorem~\ref{thm: existence of a solution of the approximated problem} with the proof in Appendix~\ref{prf: existence of a solution of the approximated problem}.

 \begin{thm}\label{thm: existence of a solution of the approximated problem}
    Consider problem \eqref{prob: approximated min covariance sensor guidance problem} and let assumptions (A5)--(A11) hold. \revision{Then problem \eqref{prob: approximated min covariance sensor guidance problem} has a solution.}
    % the first submission
    % Consider problem \eqref{prob: approximated min covariance sensor guidance problem} and let assumptions (A5)--(A11) hold. If there exists a guidance function $p_0$ in the admissible set $\mathcal{P}$ such that $\costeval{\eqref{prob: approximated min covariance sensor guidance problem}}{p_0} < \infty$, then \eqref{prob: approximated min covariance sensor guidance problem} has a solution.
\end{thm}

% \begin{pf*}{Proof}
% See Appendix~\ref{prf: existence of a solution of the approximated problem}.
% \qed
% \end{pf*}

Solving problem \eqref{prob: approximated min covariance sensor guidance problem} provides a candidate solution, denoted by $p_N^*$, where $N$ is the dimension of the approximation. The candidate $p_N^*$ may not equal the exact optimal solution, denoted by $p^*$, of the original problem~\eqref{prob: min covariance sensor guidance problem}. However, as we show in the following theorem, the candidate $p_N^*$ yields the optimal value of \eqref{prob: approximated min covariance sensor guidance problem} arbitrarily close to the one of \eqref{prob: min covariance sensor guidance problem}, as the dimension $N$ goes to infinity. Moreover, when $p_N^*$ is evaluated in the original problem \eqref{prob: min covariance sensor guidance problem}, the resulting cost is arbitrarily close to the optimal cost of \eqref{prob: min covariance sensor guidance problem}.

Before we state this convergence result, we introduce an assumption on the set of admissible guidance functions.
\begin{enumerate}
    \item[(A12)] There exist $p_{\max}>0$ and $a_{\max}>0$ such that the set of admissible guidance is $\mathcal{P}(p_{\max},a_{\max}) \defeq \{p \in C([0,t_f];P): |p(t)|$ is uniformly bounded by $p_{\max}$ and $|p(t_1) - p(t_2)| \leq a_{\max} |t_1-t_2|, \ \forall t_1,t_2 \in [0,t_f] \}$.
\end{enumerate}

Notice that the set $\mathcal{P}(p_{\max},a_{\max})$ is sequentially compact, due to the Arzel\`a-Ascoli Theorem \cite{royden2010real}, since the guidance functions in $\mathcal{P}(p_{\max},a_{\max})$ are uniformly equicontinuous and uniformly bounded.
The parameters $p_{\max}$ and $a_{\max}$ may be determined by the vehicles carrying sensors. For example, $p_{\max}$ and $a_{\max}$ refer to the maximum speed and maximum acceleration, respectively, in the case of single integrator dynamics where $p$ is the velocity command.

Theorem~\ref{thm: convergence of approximate solution minVarEst} below states the convergence of the approximation with the proof in Appendix~\ref{prf: convergence of approximate solution minVarEst}.

\begin{thm}\label{thm: convergence of approximate solution minVarEst}
    Consider problem \eqref{prob: min covariance sensor guidance problem} and its finite-dimensional approximation \eqref{prob: approximated min covariance sensor guidance problem}. Let assumptions (A4)--(A12) hold and let $p^*$ and $p_N^*$ denote the optimal guidance of \eqref{prob: min covariance sensor guidance problem} and \eqref{prob: approximated min covariance sensor guidance problem}, respectively. Then 
    \begin{equation}\label{eq: convergence of the approximate optimal cost}
        \lim_{ N \rightarrow \infty} |\optcosteval{\eqref{prob: approximated min covariance sensor guidance problem}}{p_N^*} - \optcosteval{\eqref{prob: min covariance sensor guidance problem}}{p^*}| = 0.
    \end{equation} 
    Furthermore, the cost function of \eqref{prob: min covariance sensor guidance problem} evaluated at the guidance $p_N^*$ converges to the optimal cost of \eqref{prob: min covariance sensor guidance problem}
    \begin{equation}\label{eq: convergence of the approximate optimal guidance}
        \lim_{ N \rightarrow \infty} |\costeval{\eqref{prob: min covariance sensor guidance problem}}{p_N^*} - \optcosteval{\eqref{prob: min covariance sensor guidance problem}}{p^*}| = 0.
    \end{equation} 
\end{thm}

% \begin{pf*}{Proof}
% See Appendix~\ref{prf: convergence of approximate solution minVarEst}
% \qed 
% \end{pf*}

\begin{rem}
    Two implications follow from the convergence result in Theorem~\ref{thm: convergence of approximate solution minVarEst}. First, the convergence in~\eqref{eq: convergence of the approximate optimal cost} affirms the usage of approximation since the optimal cost of the approximate problem \eqref{prob: approximated min covariance sensor guidance problem} gets arbitrarily close to that of the exact problem \eqref{prob: min covariance sensor guidance problem} as the approximation gets finer.
    Second, the convergence in \eqref{eq: convergence of the approximate optimal guidance} affirms the optimal guidance computed using approximation. When the approximate optimal guidance is evaluated by the cost function of the original problem~\eqref{prob: min covariance sensor guidance problem}, the resulting value is arbitrarily close to the optimal cost of~\eqref{prob: min covariance sensor guidance problem} as the approximation gets finer. In other words, the approximate optimal guidance is a sufficiently accurate proxy for the exact optimal guidance. %And the former is affirmed since it is numerical computable while the latter is not.
\end{rem}

\revision{
\begin{rem}
    Assumptions are made for the existence of a solution to problem \eqref{prob: min covariance sensor guidance problem} ((A9)--(A11)), well-posedness of the Riccati operators ((A1)--(A4)), and the convergence of the approximated solution ((A5)--(A8) and (A12)). Assumptions (A9)--(A11) are regarding the mobility cost and the set of admissible guidance, which are generally satisfied in engineering applications (see the discussion in \cite{Cheng2020optimalInReivew} before Theorem~3.1). The rest of the assumptions are typically satisfied with the diffusion-advection equation and Galerkin approximation (using eigenfunctions of the Laplacian operator). Details of how to check similar assumptions for the dual control problem can be found in \cite[Section~4.1]{Cheng2020optimalInReivew}.
\end{rem}

}

The convergence stated in Theorem~\ref{thm: convergence of approximate solution minVarEst} is established based on several earlier stated results, including
\begin{enumerate}
    \item the output operator's continuity with respect to location (assumption (A4)), which leads to the continuity of the uncertainty cost with respect to sensor state (Lemma~\ref{lemma: continuity of trace cost wrt sensor trajectory});
    \item existence of the Riccati operator (Lemma~\ref{lemma: existence of Riccati mild solution}) and convergence of its approximation (Theorem~\ref{thm: convergence of riccati covariance operator}); and
    \item sequential compactness of the set of admissible guidance functions (assumption (A12)), which leads to the continuity of the cost function with respect to guidance (Lemma~\ref{lem: lemma prepared for proving convergence of approximate solution} in the Appendix).
\end{enumerate}
Notice that these key results, in an analogous manner, are also required in \cite{Zhang2018sensor} when establishing the convergence to the exact optimal sensor locations of the approximate optimal locations \cite[Theorem~4.3]{Zhang2018sensor}, i.e., 
\begin{enumerate}
    \item continuity with respect to location and compactness of the output operator (dual of \cite[Theorem~2.10]{morris2010linear}), which lead to continuity of the Riccati operator with respect to sensor locations \cite[Theorem~4.1]{Zhang2018sensor};
    \item existence of the Riccati operator \cite[Theorem~2.5]{Zhang2018sensor} and the convergence of its approximation \cite[Theorem~4.2]{Zhang2018sensor}; and
    \item sequential compactness of the set of admissible locations, which is inherited from the setting that the spatial domain is closed and bounded in a finite-dimensional space.
\end{enumerate}

% {\color{red} I think we should focus more on the technical difference than the details. Say ... are difference rather than the obvious finite-dimensional vs. infinite-dimensional or the mobile vs. stationary crap.}

% Although the establishment of convergence is similar to the one in \cite{morris2010linear}, the cost function and type of Riccati equation are different: we have quadratic PDE cost plus generic mobility cost and differential Riccati equation in this paper for control and actuator guidance versus Riccati operator's norm as cost function and algebraic Riccati equation in \cite{morris2010linear} for actuator placement.

To compute an optimal solution to problem \eqref{prob: approximated min covariance sensor guidance problem}, we use Pontryagin's minimum principle, which introduces a Hamiltonian function that has the same form as \eqref{eq: Hamiltonian} except that the covariance operator $\Pi(t)$ in \eqref{eq: Hamiltonian} is replaced by its approximation $\Pi_N(t)$. And correspondingly, the resulting two-point boundary value problem has the same form as \eqref{eq: pontryagin's principle} except for $\Pi(t)$ being replaced by $\Pi_N(t)$.

\revision{
\begin{rem}
    The optimal sensor trajectory $\zeta^*$ steered by the optimal guidance $p^*$ may be used as the reference trajectory tracked by the vehicle’s lower-level control. Although collision avoidance among the sensors is not discussed in this paper, it can be incorporated into the lower-level control using numerous methods in the existing literature, e.g., \cite{wang2017safety} and the references therein.
\end{rem}
}

\section{Simulation results}
This section shows the simulation results obtained using the solution method proposed in Section~\ref{sec: approximation}. Comparison and analysis are made regarding the performance of the mobile sensor(s) under optimal guidance for the case of a single sensor, a team of homogeneous sensors, and a team of heterogeneous sensors.

\begin{figure}[t]
	\centering
% 	\vspace{-0.2cm}
	\includegraphics[width=\columnwidth]{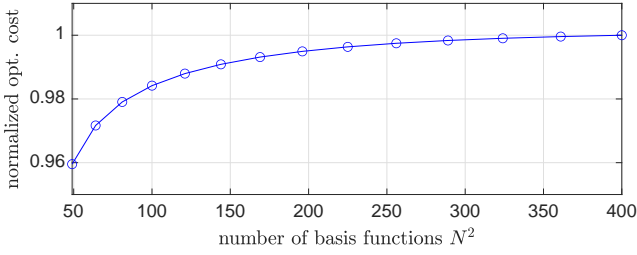}
	\caption{Approximate optimal costs $\optcosteval{\eqref{prob: approximated min covariance sensor guidance problem}}{p_N^*}$ normalized with respect to the optimal cost for $N^2 = 400$.}
	\label{fig: convergence}
\end{figure}

We use the Galerkin scheme to approximate the infinite-dimensional terms. The orthonormal set of eigenfunctions of the Laplacian operator $\nabla^2$ (with zero Dirichlet boundary condition) over the spatial domain $\Omega = [0,1] \times [0,1]$ is $\phi_{i,j}(x,y) = 2 \sin(\pi i x) \sin(\pi j y)$. With a single index $k \defeq (i-1)N + j$ such that $\phi_k \defeq \phi_{i,j}$, the set of eigenfunctions $\{\phi_k \}_{k=1}^{N^2}$ spans an $N^2$-dimensional space that is previously denoted by $\mathcal{H}_N$. For the orthogonal projection $P_N: \mathcal{H} \rightarrow \mathcal{H}_N$, it follows that $P_N^{\star} = P_N$ and $P_N^{\star} P_N \rightarrow I$ strongly. % \cite{Burns2015solutions} 
And the assumption~\eqref{eq: primal and dual convergence of the approximation} in Theorem~\ref{thm: convergence of riccati covariance operator} holds uniformly for all $t \in [0,t_f]$ as $N \rightarrow \infty$~\cite{Burns2015solutions}. \revision{We set $p_{\max}$ and $a_{\max}$ to be sufficiently large so that the solution is in the set $\mathcal{P}(p_{\max},a_{\max})$. We plot the optimal cost $\optcosteval{\eqref{prob: approximated min covariance sensor guidance problem}}{p_N^*}$ for $N$ from 7 to 20, as shown in Fig.~\ref{fig: convergence}. The optimal cost shows a tendency of exponential convergence as we increase the number of basis functions. And we choose $N=12$ in the rest of the simulations since it is the smallest dimension with the optimal cost within 1\% of the optimal cost evaluated with the maximum dimension $N = 20$ in the trials.}

The parameters in the simulation are $t_f = 2$, and $a=0.01$.
We use single integrator dynamics for each sensor. The state $\zeta$ is the 2D location of the sensors and guidance $p$ is the 2D velocity command. In some applications, the flow field $\mathbf{v}$ of the diffusion-advection process can affect the mobile sensors. For example, surface vehicles that measure the concentration of certain chemical substances or biological entities in a water body are subject to the movement of the water. Considering this realistic condition, we append the flow field $\mathbf{v} = [0.1,-0.1]^\top$ of advection to the right-hand side of the single integrator dynamics, which means the sensors will drift along the flow when zero guidance is implemented. The previous statements and results on the existence of solution and convergence of the approximate solution still hold within this setting. The optimization will find optimal guidance subject to (or possibly taking advantage of) this flow field. The sensor has the square-shaped average output kernel (see \eqref{eq: examplery measurement kernel function}) with $r_i= 0.05$, in which case its footprint covers only $1\%$ of the domain in area.

We set $g(p(t)) = \gamma p^\top(t) p(t)/2$ and $h(\zeta(t),t) = h_f(\zeta(t_f)) = 0$ as the mobility cost, which is simply the quadratic guidance effort for $\gamma > 0$.

Assumption (A4) holds for the choice of output operator (see Remark~\ref{remark: continuity wrt location}). With the Galerkin approximation using the orthonormal eigenfunctions $\{\phi_k \}_{k=1}^{N^2}$, it can be shown that assumption (A8) holds for $l_N (\cdot) = N^2 l(\cdot)$. 
Assumptions (A5)--(A7) hold with $q=1$ under the Galerkin approximation with aforementioned basis functions $\{ \phi_k\}_{k=1}^{N^2}$ \cite{Burns2015solutions}. Assumptions (A9)--(A11) and (A12) hold for the choice of functions in the mobility cost and parameters of the set of admissible guidance functions, respectively.

To evaluate the performance of the optimal guidance in simulation, we set the deterministic portion of the initial condition of the PDE to zero, i.e., $\hat{Z}_0 = 0$, which excludes the bias from choosing a particular non-zero one. The stochastic portion of the initial condition $w_0$ is chosen as a zero-mean Gaussian process with non-stationary kernel function $k_0: \Omega \times \Omega \rightarrow \mathbb{R}^+$ such that $k_0(x_1,x_2) \defeq 9 \text{exp}(-|x_{12}|_2^2/200- |x_{10}|_2^2/10 - |x_{20}|_2^2/10)$,
% \begin{equation}\label{eq: initial uncertainty kernel}
%     k_0(x_1,x_2) \defeq 9 \text{exp}\Big(-\frac{|x_{12}|_2^2}{200}- \frac{|x_{10}|_2^2}{10} - \frac{|x_{20}|_2^2}{10} \Big),
% \end{equation}
where $x_{ij}\defeq x_i-x_j$ for $i,j \in \{0,1,2\}$ and $x_0 \in \Omega$ represents the peak of the uncertainty in the domain and we set $x_0 = [0.75,0.25]^\top$. For the state noise $w(t)$, we use a zero-mean Gaussian process with a homogeneous kernel function $k: \Omega \times \Omega \rightarrow \mathbb{R}^+$ such that $k(x_1,x_2) = \text{exp}(-|x_{12}|_2^2/2000 )$.

We use the forward-backward sweep method \cite{mcasey2012convergence} to solve the two-point boundary problem \eqref{eq: pontryagin's principle} (with $\Pi$ replaced by $\Pi_N$) and subsequently compute the optimal guidance. A fixed-step length of 0.01 and a relative tolerance of $1 \times 10^{-6}$ are applied in the iterative procedure.
The forward propagation of \eqref{eq: Pontryagin vehicle dynamics} and \eqref{eq: mild solution of approximate Riccati covariance} and the backward propagation of \eqref{eq: Pontryagin costate dynamics} are computed via the Runge-Kutta method.

\subsection{Single sensor results}\label{sec: single sensor simulation}
Two important parameters in the problem setting are the sensor noise variance $R$ and mobility penalty $\gamma$. Smaller $R$ yields higher sensor quality, whereas smaller $\gamma$ yields better mobility of the vehicle. For example, if $\gamma$ is the mass of the vehicle, then the guidance effort is the kinetic energy of the vehicle. These parameters affect the performance of the estimation as shown next. %The terminologies \textit{uncertainty cost} and \textit{guidance effort} refer to $\int_0^{t_f} \trace{\Pi(t)} \dd t$ and $\frac{1}{2}\int_0^{t_f} p(t)^\top \gamma p(t) \dd t$, respectively.
In this simulation, the sensor is initiated at $\zeta_0 = [0.3,0.1]^\top$. Monte Carlo simulations of 100 trials compare the optimal guidance with three naive guidance policies as follows: naive~1 crosses the domain by reaching the opposite of the initial location within domain at $[0.7 ,0.9]^\top$ at a constant speed; naive~2 reaches the peak of the initial uncertainty $x_0$ at a constant speed; and naive~3 circulates the domain in the clockwise direction with center $[0.5,0.5]^\top$, radius $1/\sqrt{5}$, and angular speed $\pi$ rad/s.
% \begin{enumerate}
%     \item naive 1: crossing the domain by reaching the opposite of the initial location within domain at $[0.7 ,0.9]^\top$ at a constant speed.
%     \item naive 2: reaching the peak of the initial uncertainty $x_0$ at a constant speed.
%     \item naive 3: circulating the domain in the clockwise direction with center $[0.5,0.5]^\top$, radius $1/\sqrt{5}$, and angular speed $\pi$ rad/s.
% \end{enumerate}
A stationary sensor is also included for comparison whose null guidance merely compensates for the flow field $\mathbf{v}$. 
\begin{figure}[t]
	\centering
	\begin{subfigure}[b]{\columnwidth}
		\includegraphics[width=\textwidth]{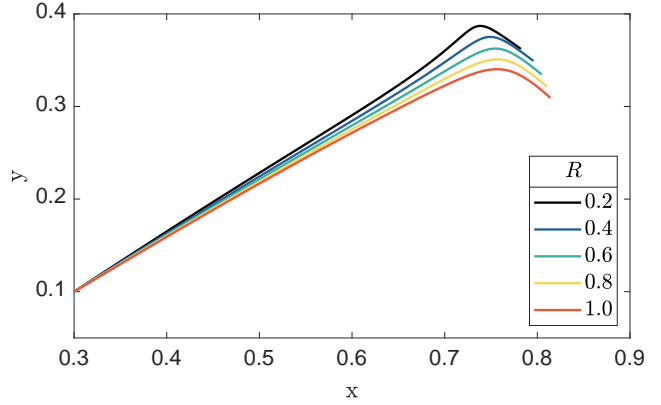}
		\subcaption{Optimal trajectories for various values of the sensor noise variance $R$}
		\label{fig: traj_fixed_gamma}
	\end{subfigure}\\
% 	\begin{subfigure}[b]{\columnwidth}
% 		\includegraphics[width=\textwidth]{costBreakdown_fixed_gamma.eps}
% 		\subcaption{Optimal cost and its breakdown}
% 		\label{fig: costBreakdown_fixed_gamma}
% 	\end{subfigure}\\
	\begin{subfigure}[b]{\columnwidth}
		\includegraphics[width=\textwidth]{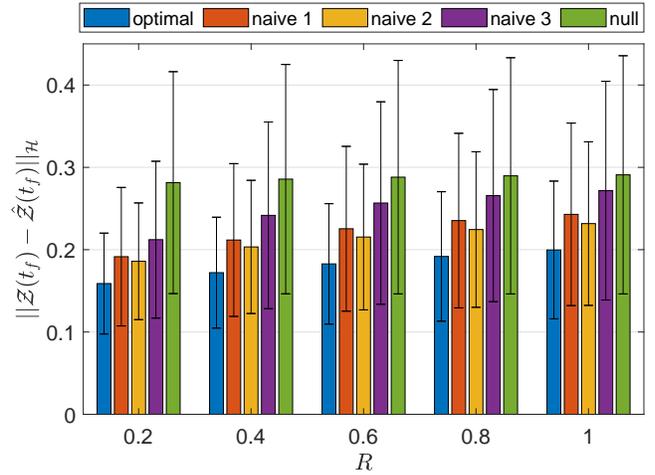}
		\subcaption{Norm of the terminal estimation error for various values of sensor noise variance $R$. The color bar shows the mean value; the error bar shows the standard deviation.}
		\label{fig: MC_fixed_gamma}
	\end{subfigure}
	\caption{The value of mobility penalty $\gamma$ is fixed at 0.5, whereas the sensor noise's variance $R$ takes values in the set $\{0.2,0.4,0.6,0.8,1 \}$.}
\end{figure}

First, hold either $R$ or $\gamma$ fixed and vary the other to observe the variation of the optimal trajectory. 
Fig. \ref{fig: traj_fixed_gamma} displays the trajectories when $\gamma = 0.5$ and $R$ varies from $0.2$ to $1$. The sensor maneuvers less as $R$ increases, which indicates the optimal guidance's compensation for deteriorating sensor quality by moving it closer to the peak of the initial uncertainty at $x_0 = [0.75,0.25]^\top$.

For the Monte Carlo trials, the mean and standard deviation of the terminal estimation error's norm are shown in Fig. \ref{fig: MC_fixed_gamma}. The optimal guidance exhibits smaller mean and variance of terminal estimation error over the naive guidance policies and the null guidance at each evaluated $R$. Notice that the advantage is preserved when the sensor quality deteriorates as $R$ increases.

\begin{figure}[t]
	\centering
	\begin{subfigure}[b]{\columnwidth}
		\includegraphics[width=\textwidth]{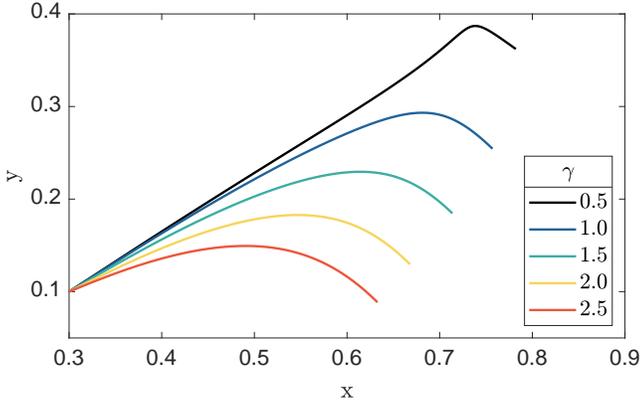}
		\subcaption{Optimal trajectories for various values of the mobility penalty $\gamma$}
		\label{fig: traj_fixed_R}
	\end{subfigure}\\
% 	\begin{subfigure}[b]{\columnwidth}
% 		\includegraphics[width=\textwidth]{costBreakdown_fixed_R.eps}
% 		\subcaption{Optimal cost and its breakdown}
% 		\label{fig: costBreakdown_fixed_R}
% 	\end{subfigure}\\
	\begin{subfigure}[b]{\columnwidth}
		\includegraphics[width=\textwidth]{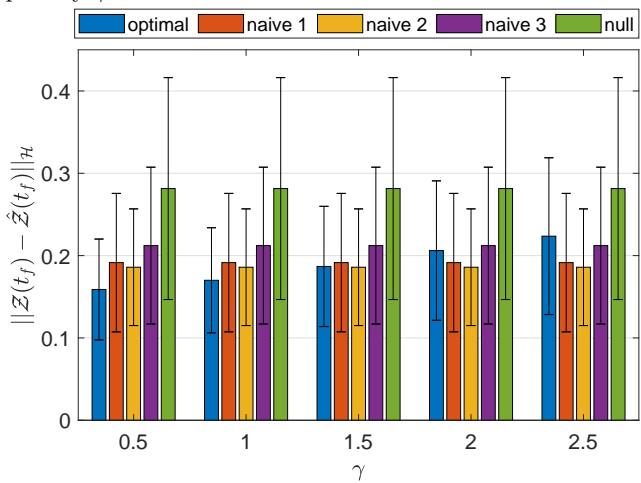}
		\subcaption{Norm of the terminal estimation error for various values of the mobility penalty $\gamma$. The color bar shows the mean value; the error bar shows the standard deviation.}
		\label{fig: MC_fixed_R}
	\end{subfigure}%\\
	\caption{The value of sensor noise's variance $R$ is fixed at $0.2$, whereas the mobility penalty $\gamma$ takes values in the set $\{0.5,1,1.5,2,2.5 \}$.}
	\label{fig: fixed_R}

\end{figure}

Fig.~\ref{fig: traj_fixed_R} displays the trajectories when $R = 0.2$ and $\gamma$ varies from $0.5$ to $2.5$. A bigger value of $\gamma$ suppresses the guidance effort and hence reduces the range of the sensor's motion. The results of Monte Carlo trials when $\gamma$ varies are shown in Fig.~\ref{fig: MC_fixed_R}. Notice that the mean and standard deviation are invariant for each naive guidance policy and the null guidance among the varying $\gamma$ since the sensor trajectory steered by each of these guidance policies is independent of $\gamma$. The optimal guidance shows its advantage over the other guidance policies at a relatively smaller values of $\gamma$, e.g., at $0.5$ and $1$. This advantage is gradually lost as $\gamma$ takes relatively bigger values, e.g., $2$ and $2.5$. This comparison suggests that the optimal guidance may not be the best option for vehicles with large $\gamma$ (e.g., when the vehicle is heavy), despite the fact that the guidance is still optimal for the chosen cost function in simulation.

\begin{figure*}[t]
	\centering
	\includegraphics[width=\textwidth]{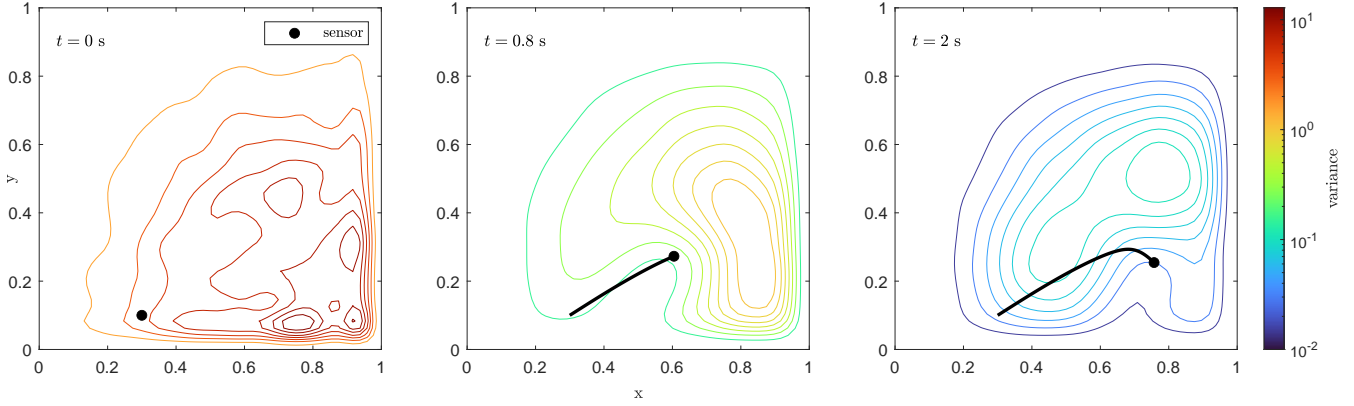}
	\caption{Trajectory of the sensor following the optimal guidance. The contours represent the level curves of pointwise variance of the estimation error on a uniform grid within the domain, computed from the Monte Carlo simulations.}
	\label{pic: uncertainty evolution}
\end{figure*}

Fig.~\ref{pic: uncertainty evolution} shows the snapshots of the sensor trajectory under the optimal guidance and contour plot of the pointwise variance of the estimation error among the Monte Carlo trials with $\gamma = 1$ and $R = 0.2$. The pointwise variance is computed at each point in a uniform grid of $144 \times 144$ sampling points in the domain $\Omega$. The sensor is steered quickly towards the area with higher uncertainty near $[0.75, 0.25]^\top$ (see the snapshot at $t = 0.8$ s). Eventually, the sensor starts to drift along the flow field $\mathbf{v} = [0.1,-0.1]^\top$ (see the snapshot at $t = 2$ s). That sensor can effectively reduce the uncertainty of the estimation error as can be observed from the drop in the pointwise variance in the sensor's footprint. Note that the zero Dirichlet boundary condition also contributes to reducing the uncertainty of the estimation error via diffusion and advection.

\subsection{Team of homogeneous sensors}

\begin{figure*}[t]
	\centering
	\includegraphics[width=\textwidth]{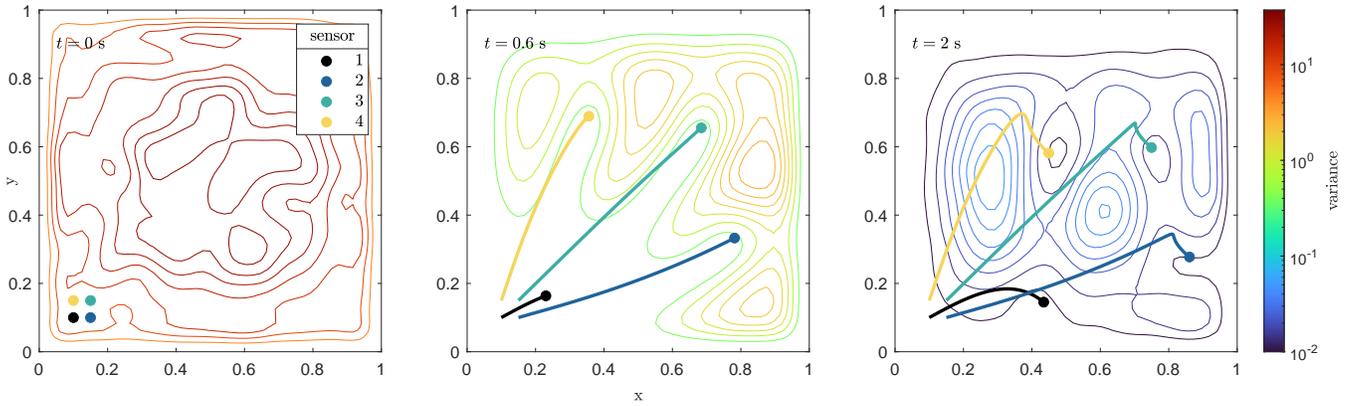}
	\caption{Trajectories of multiple homogeneous sensors following the optimal guidance. The contours represent the level curves of pointwise variance of the estimation error on a uniform grid within the domain. }
	\label{pic: uncertainty evolution 4 agents}
\end{figure*}
To demonstrate the framework's capability in guiding a team of multiple sensors, we simulate four homogeneous sensors ($R = 0.2$, $\gamma = 0.5$). To adapt to a total of four sensors, the kernel functions for $w_0$ and $w$ are set to be $4k_0(x_1,x_2)$ and $4k(x_1,x_2)$, respectively, 
and the peak of the initial uncertainty is set to the center $[0.5,0.5]^\top$. The other settings are identical to those in Section~\ref{sec: single sensor simulation}. Fig.~\ref{pic: uncertainty evolution 4 agents} shows the snapshots of the sensors' trajectories under the optimal guidance and contour plot of the pointwise variance of the estimation error among the Monte Carlo trials. Similar to the case of a single sensor in Fig.~\ref{pic: uncertainty evolution}, the sensors quickly sweep the peak of the initial uncertainty in the center and expand to cover the domain (see the snapshot at $t = 0.6$ s). The pointwise variance drops along sensors' footprints. The sensors essentially drift along the flow to reduce the guidance effort (see the snapshot at $t = 2$ s).

\subsection{Team of heterogeneous sensors}

The parameters $R$ and $\gamma$ essentially relate to operational planning: one may invest more for better sensor quality or a swifter vehicle. Consequently, one would necessarily invest more for a team of superior mobile sensors (e.g., $R = 0.2$ and $\gamma = 0.5$) than a team of poor mobile sensors (e.g., $R = 1$ and $\gamma = 2.5$). The latter has five times as much sensor noise (in terms of standard deviation) and five times the mobility penalty as the former. One may balance the conflicting needs of performance and investment by deploying a team of heterogeneous sensors, i.e., a mixed team of superior and poor sensors. 

\begin{figure}[t]
	\centering
% 	\vspace{-0.2cm}
	\includegraphics[width=\columnwidth]{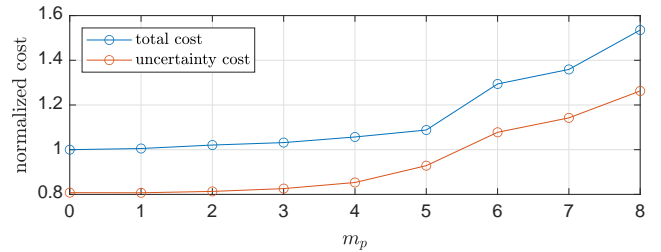}
	\caption{Normalized optimal total cost and uncertainty cost of a heterogeneous team with $m_p$ poor sensors and $8-m_p$ superior sensors.}% The dashed lines indicate a $20\%$ increment to the optimal total cost (blue) and uncertainty cost (red) of the homogeneous team of superior sensors ($m_p = 0$).}
	\label{fig: one best plus m_w worst}
\end{figure}

The following simulation compares the performance of a team of heterogeneous sensors (including $m_p$ poor sensors and $8-m_p$ superior sensors, for $m_p$ in the range of $\{1,2,\dots,7\}$) with that of homogeneous teams ($m_p = 0$ for superior sensors only and $m_p = 8$ for poor sensors only). The sensors are initiated in the lower left corner of the domain. % and the coordinates are shown in Table~\ref{tb: the table for initial location}. 
To adapt to a total of eight sensors, the kernel functions for $w_0$ and $w$ are set to be $8k_0(x_1,x_2)$ and $8k(x_1,x_2)$, respectively, and the peak of the initial uncertainty is set to the center of the domain at $[0.5,0.5]^\top$. Fig.~\ref{fig: one best plus m_w worst} shows the normalized optimal total cost and uncertainty cost for the heterogeneous team compared with homogeneous teams of superior and poor sensors. The performance deteriorates judged by the rising costs as the number of poor sensors increases in team. However, the degradation of the heterogeneous team (when $m_p \leq 5$) is maintained within $20\%$ of the superior team in both the total cost and uncertainty cost, which indicates the cost effectiveness of the heterogeneous team since the investment reduces linearly as $m_p$ increases.

% \renewcommand{\arraystretch}{1} % Default value: 1
%   \captionsetup{%size=footnotesize,
% 	%justification=centering, %% not needed
% 	skip=5pt, position = bottom}
% \begin{table}[h]
%     \centering
% 	\small
% 	%\captionsetup{font=small}
% 	\caption{Coordinates of initial sensor locations}
% 	\begin{tabular*}{\columnwidth}{@{\extracolsep{\fill}}ccccccccc} 
% 		\toprule[1pt]
% 		sensor & 1 & 2 & 3 & 4 & 5 & 6 & 7 & 8 \\
% 		\midrule
% 		x  & 0.1 & 0.15 & 0.2 & 0.1 & 0.15 & 0.2 &  0.1 & 0.15 \\
%         y  & 0.1 & 0.1 & 0.1 & 0.15 &  0.15 & 0.15 &  0.2 & 0.2\\
% 		\bottomrule[1pt]
% 	\end{tabular*}\label{tb: the table for initial location}
% % 	\vspace{-0.3cm}
% \end{table}
% \normalsize

\section{Conclusion}
This paper proposes a guidance design method for a team of mobile sensors to estimate a spatiotemporal process modeled by a 2D diffusion-advection process. We formulate an optimization problem that minimizes the sum of the trace of the covariance operator of the Kalman-Bucy filter and a generic mobility cost of the sensors subject to the dynamics of the sensor platforms. Conditions for the existence of a solution to the proposed problem are established, where a fundamental assumption is that the output kernel is continuous with respect to location. Approximation of the infinite-dimensional terms permits the computation of the optimal guidance, and we prove that the approximated problem's optimal cost converges to that of the exact problem. Moreover, we prove that the optimal guidance obtained from an approximate problem yields the cost evaluated by the exact problem's cost function arbitrarily close to the exact optimal cost. %, which affirms the usage of approximation in obtaining optimal guidance.
Exponential convergence in the optimal cost can be observed in simulations when the approximation dimension increases.

We use Pontryagin's minimum principle to numerically compute optimal guidance. The numerical solutions are evaluated in simulations. We study how the optimal trajectory and terminal estimation error change subject to the varying values of sensor noise variance and the mobility penalty for a single sensor. Trajectories of the sensors and the evolution of the pointwise variance are shown for a single sensor and a homogeneous team of sensors in a flow field. We also study the cost-effectiveness of a heterogeneous team of mixed sensors with superior and poor qualities by comparing the performance degradation to that of the homogeneous team of superior sensors.
%Simulation results show that, for a single sensor, either improving sensor quality or reducing mobility penalty improves the performance evaluated by the proposed cost function and terminal estimation error. For a team of mobile sensors, degradation of the estimation performance is kept within 20\% when several superior sensors are replaced by the poor alternatives, indicating the cost effectiveness of a heterogeneous team.

Ongoing and future work includes establishing a convergence rate in Theorem~\ref{thm: convergence of approximate solution minVarEst}, incorporating constraints to the problem formulation (e.g., constraints on sensor states), \revision{extension to more complicated sensor dynamics (e.g., control-affine),} and demonstrating the framework experimentally using a swarm of quadrotor helicopters in an outdoor testbed.

\section*{Acknowledgments} 
The authors would like to thank the reviewers for their helpful suggestions. Work reported in this article was sponsored under a seed grant from Northrop Grumman.

\appendix

\section{Proof of Theorem~\ref{thm: existence of a solution of minVarGuidanceProblem}}\label{prf: existence of a solution of minVarGuidanceProblem}
\vspace{-0.3cm}

Since we have proved that the uncertainty cost $\int_0^{t_f} \trace{\Pi(t)} \dd t$ is a continuous mapping $K(\cdot)$ of the sensor state $\zeta$ (see Lemma~\ref{lemma: continuity of trace cost wrt sensor trajectory}), the existence of optimal guidance of \eqref{prob: min covariance sensor guidance problem} can be proved using the techniques of proving existence of solution to an optimal control (guidance) problem. %for which we refer to Definition~\ref{def: applied for existence of optimization solution} and Theorem~\ref{thm: optimization existence of solution} stated below.

% \begin{defn}\label{def: applied for existence of optimization solution}
%     Suppose $(X,\norm{\cdot})$ is a normed linear space. 
%     \begin{enumerate}
%         \item A sequence $\{x_k \} \subset X$ is weakly convergent to $x \in X$, denoted by $x_k \rightharpoonup x$, if $\lim \limits_{k \rightarrow \infty} (x^{\star},x_k) = (x^{\star},x)$ for all $x^{\star}$ belonging to the dual space $X^{\star}$. 
%         \item A subset $U \subset X$ is weakly sequentially closed if $\{x_k \} \subset U$ and $x_k \rightharpoonup x$ implies $x \in U$.
%         \item A subset $U \subset X$ is weakly sequentially compact if for every sequence $\{x_k \} \subset U$ there exists a subsequence $\{x_{k_i} \} \subset \{x_k \}$ and an $x \in U$ with $x_{k_i} \rightharpoonup x$.
%         \item Suppose $U \subset X$ and $f:U \rightarrow \mathbb{R}$. The mapping $f$ is weakly sequentially lower semicontinuous on $U$ if $\{ x_k\} \subset U$ and $x_k \rightharpoonup x \in U$ implies $f(x) \leq \liminf_{k \rightarrow \infty} f(x_k)$.
%     \end{enumerate}
% \end{defn}

% \begin{thm}[\protect{\cite[Theorem 6.1.4]{werner2013optimization}}]\label{thm: optimization existence of solution}
%     Let $(X,\norm{\cdot})$~be a normed linear space, $U_0 \subset X$ be weakly sequentially compact and $f : U_0 \rightarrow \mathbb{R}$ be weakly sequentially lower semicontinuous on $U_0$. Then there exists an $\bar{x} \in U_0$ such that
%     \begin{equation}
%         f(\bar{x}) = \inf \{f(x): x \in U_0 \}.
%     \end{equation}
% \end{thm}

\begin{pf*}{Proof of Theorem~\ref{thm: existence of a solution of minVarGuidanceProblem}}
%The proof is based on \cite[p.~219]{werner2013optimization}. 
Without loss of generality, we consider the case of one mobile sensor, i.e., $m_s = 1$. The case of $m_s \geq 2$ follows naturally.

Our proof follows the proof in \cite[Chapter~6.2]{werner2013optimization} which proves existence of a solution to an optimal control problem based on functional analytic theorems. Consider problem \eqref{prob: min covariance sensor guidance problem}'s admissible set of guidance functions $\mathcal{P} \defeq \{p \in L^2([0,t_f];\mathbb{R}^m): p(t) \in P, t \in [0,t_f] \}$. 
\revision{Since there exists $p_0 \in \mathcal{P}$ such that $\costeval{\eqref{prob: min covariance sensor guidance problem}}{p_0} < \infty$ (e.g., $p_0=0$ that yields a stationary sensor at $\zeta_0$)}, let $\mathcal{P}_0 \defeq \{p \in \mathcal{P}: \ \costeval{\eqref{prob: min covariance sensor guidance problem}}{p} \leq \costeval{\eqref{prob: min covariance sensor guidance problem}}{p_0} \}$.
% in the first submission:
% Since we have assumed the existence of $p_0 \in \mathcal{P}$ such that $\costeval{\eqref{prob: min covariance sensor guidance problem}}{p_0} < \infty$ (recall that the notation $\costeval{\eqref{prob: min covariance sensor guidance problem}}{p_0}$ stands for the cost function of \eqref{prob: min covariance sensor guidance problem} evaluated at $p_0$), let $\mathcal{P}_0 \defeq \{p \in \mathcal{P}: \ \costeval{\eqref{prob: min covariance sensor guidance problem}}{p} \leq \costeval{\eqref{prob: min covariance sensor guidance problem}}{p_0} \}$. We wish to prove \underline{Condition-1}, \underline{Condition-2}, and \underline{Condition-3} stated below:

    % \vspace{0.1cm}
    \noindent \underline{Condition-1}: The set $\mathcal{P}_0$ is bounded.
    
    \noindent \underline{Condition-2}: The set $\mathcal{P}_0$ is weakly sequentially closed.
    
    \noindent \underline{Condition-3}: The mapping $\costeval{\eqref{prob: min covariance sensor guidance problem}}{\cdot}: \mathcal{P} \rightarrow \mathbb{R}^+$ is weakly sequentially lower semicontinuous on $\mathcal{P}_0$.
    
% And by Theorem~\ref{thm: optimization existence of solution}, if \underline{Condition-1} and \underline{Condition-2} hold then there exists a guidance function in $\mathcal{P}_0$ such that the cost function of \eqref{prob: min covariance sensor guidance problem} is minimized.
\underline{Condition-1} and \underline{Condition-2} imply that $\mathcal{P}_0$ is weakly sequentially compact. By \cite[Theorem~6.1.4]{werner2013optimization}, problem \eqref{prob: min covariance sensor guidance problem} has a solution when \underline{Condition-1}--\underline{Condition-3} hold.

Before showing these three conditions, we first define a the solution map $T: L^2([0,t_f];\mathbb{R}^m) \rightarrow C ([0,t_f];\mathbb{R}^n)$ of the sensor dynamics \eqref{eq: general dynamics of the mobile sensor} by $(Tp)(t) \defeq \zeta(t) = e^{\alpha t} \zeta_0 + \int_0^t e^{\alpha (t-\tau)} \beta p(\tau) \dd \tau $ for $t \in [0,t_f]$. 
\revision{The continuity of the map $T$ is straightforward \cite{werner2013optimization}, i.e., there exists $c_4>0$ such that $\norm{Tp}_{C([0,t_f];\mathbb{R}^n)} \leq c_4 \norm{p}_{L^2([0,t_f];\mathbb{R}^m)}$.
}
% For $p_1,p_2 \in L^2([0,t_f];\mathbb{R}^m)$ and $t \in [0,t_f]$, we have 
% \begin{align}
%     & |Tp_1(t) - Tp_2(t)|_1 \nonumber \\
%      \leq & \int_0^{t} |e^{\alpha (t-\tau)} \beta|_1 |p_1(\tau) - p_2(\tau)|_1 \dd \tau \nonumber\\
%     \leq & c_4 \int_0^{t} |e^{\alpha (t-\tau)} \beta|_1 |p_1(\tau) - p_2(\tau)|_2 \dd \tau \nonumber\\
%     % \leq & c_4 \left( \int_0^{t} |e^{\alpha (t-\tau)} \beta|_1^2 \dd \tau \right)^{1/2} \left(\int_0^t |p_1(\tau) - p_2(\tau)|_2^2 \dd \tau \right)^{1/2} \nonumber \\
%     = & c_4 \left( \int_0^{t} |e^{\alpha (t-\tau)} \beta|_1^2 \dd \tau \right)^{1/2} \norm{p_1-p_2}_{L^2([0,t];\mathbb{R}^m)} \nonumber \\
%     \leq & c_4 c_5 \norm{p_1-p_2}_{L^2([0,t];\mathbb{R}^m)}
% \end{align}
% for  $c_4$ and $c_5 > 0$. Hence, 
% \begin{align}
%      \norm{Tp_1 - Tp_2}_{C([0,t_f];\mathbb{R}^n)} = & \sup_{t \in [0,t_f]} |Tp_1(t) - Tp_2(t)|_1 \nonumber \\
%      \leq & c_4 c_5 \norm{p_1-p_2}_{L^2([0,t_f];\mathbb{R}^m)},
% \end{align}
% which also shows that $T$ is a continuous mapping, i.e., for all $p \in L^2([0,t_f];\mathbb{R}^m)$
% \begin{equation}\label{eq: boundedness of the mapping T}
%     \norm{Tp}_{C([0,t_f];\mathbb{R}^n)} \leq c_4 c_5 \norm{p}_{L^2([0,t_f];\mathbb{R}^m)}.
% \end{equation}

% \vspace{0.1cm}
\noindent \textit{\textbf{Proof of \underline{Condition-1}}}: Suppose $p \in \mathcal{P}_0$, then
\begin{align}
    \costeval{\eqref{prob: min covariance sensor guidance problem}}{p_0} \geq & \ \costeval{\eqref{prob: min covariance sensor guidance problem}}{p} \nonumber \\
    = & \textstyle \int_0^{t_f} h(Tp(t),t) + g(p(t),t)  + \trace{\Pi(t)} \dd t    \nonumber \\
    & + h_f(Tp(t_f)) \nonumber \\
    \geq & \ \textstyle \int_0^{t_f} d_1 |p(t)|^2_2 \dd t \nonumber \\
    = & \ d_1 \norm{p}_{L^2([0,t_f];\mathbb{R}^m)}^2.
\end{align}
Since $d_1>0$, the boundedness of $\mathcal{P}_0$ follows.
% First, we define the norm on $L^2([0,t];\mathbb{R}^m)$: for $p \in L^2([0,t];\mathbb{R}^m)$,
% \begin{equation}
%     \norm{p}_{L^2([0,t];\mathbb{R}^m)} =  \left( \int_0^{t} |p(\tau)|^2_2 \dd \tau \right)^{1/2}.
% \end{equation}
% The space $L^2([0,t];\mathbb{R}^m)$ equipped with norm $\norm{\cdot}_{L^2([0,t];\mathbb{R}^m)}$ is a reflexive Banach space for $t \in [0,t_f]$ \cite{yosida1988functional}. % see Example 3 on p. 115.

% To show that the set $\mathcal{P}_0$ is weakly sequentially compact, by \cite[Theorem~2.11]{troltzsch2010optimal}, %see google books:
% % Theorem~2.11: Every convex and closed subset of a Banach space is weakly sequentially closed. If the space is reflexive and the set is in addition bounded, then it is weakly sequentially compact.
% it suffices to show that the set $\mathcal{P}_0$ is bounded because $\mathcal{P}_0$ is convex and closed (both inherited from the set $P$) and is a subset of a reflexive Banach space $L^2([0,t_f];\mathbb{R}^m)$.

\noindent \textit{\textbf{Proof of \underline{Condition-2}}}: Suppose $\{p_k \} \subset \mathcal{P}_0$ and $\{p_k \}$ converges to $p$ weakly (denoted by $p_k \rightharpoonup p$). We wish to show $p \in \mathcal{P}_0$. We start with establishing that $\mathcal{P}$ is weakly sequentially closed and, hence, $p \in \mathcal{P}$. Subsequently, we show $\costeval{\eqref{prob: min covariance sensor guidance problem}}{p} \leq \costeval{\eqref{prob: min covariance sensor guidance problem}}{p_0}$ to conclude \underline{Condition-2}.

To show that the set $\mathcal{P}_0$ is weakly sequentially closed, by \cite[Theorem~2.11]{troltzsch2010optimal}, %see google books:
% Theorem~2.11: Every convex and closed subset of a Banach space is weakly sequentially closed. If the space is reflexive and the set is in addition bounded, then it is weakly sequentially compact.
it suffices to show that $\mathcal{P}$ is closed and convex. Let $\{q_k\} \subset \mathcal{P}$ and $q_k \rightarrow q$. We want to show $q \in \mathcal{P}$, i.e., $q \in L^2([0,t_f];\mathbb{R}^m)$ and $q(t) \in P$ for $t \in [0,t_f]$. Since $L^2([0,t_f];\mathbb{R}^m)$ is complete, we can choose a subsequence $\{q_{k_j} \} \subset \mathcal{P}$ that converges to $q$ pointwise almost everywhere on $[0,t_f]$ \cite[p.~53]{yosida1988functional}. Since $P$ is closed (assumption (A9)), $q(t) \in P$ for almost all $t \in [0,t_f]$. Hence, $\mathcal{P}$ is closed. The convexity of $\mathcal{P}$ follows from that of $P$ (assumption (A9)), i.e., if $p_1,p_2 \in \mathcal{P}$, then $\lambda p_1 + (1- \lambda)p_2 \in L^2([0,t_f];\mathbb{R}^m)$ and $\lambda p_1(t) + (1- \lambda)p_2(t) \in P$ for $t \in [0,t_f]$ and $\lambda \in [0,1]$.

What remain to be shown is $\costeval{\eqref{prob: min covariance sensor guidance problem}}{p} \leq \costeval{\eqref{prob: min covariance sensor guidance problem}}{p_0}$. Since $p_k \rightharpoonup p$, by definition, we have $Tp_k \rightarrow Tp$. We now show that the sequence $\{Tp_k\}$ contains a uniformly convergent subsequence in $C([0,t_f];\mathbb{R}^n)$. 
The sequence $\{Tp_k \} \subset C([0,t_f];\mathbb{R}^n)$ is uniformly bounded and uniformly equicontinuous for the following reasons: Since $\norm{Tp_k}_{C([0,t_f];\mathbb{R}^n)} \leq c_4  \norm{p_k}_{L^2([0,t_f];\mathbb{R}^m)}$, it follows that $\norm{Tp_k}_{C([0,t_f];\mathbb{R}^n)}$ is uniformly bounded, because $\{p_k\} \subset \mathcal{P}_0$ which is a bounded set. For $s,t \in [0,t_f]$, we have%there exists $c_7>0$ and $c_8>0$ such that 
\begin{align*}
     |Tp_k(s) - Tp_k(t)|_1  = &\ \left|\textstyle \int_s^t \alpha Tp_k(\tau) +  \beta  p_k(\tau) \dd \tau \right|_1 \\
    \leq &\ |t-s||\alpha|_1\norm{Tp_k}_{C([0,t_f];\mathbb{R}^n)} \nonumber \\
    & + |t-s|^{1/2} |\beta|_2 \norm{p_k}_{L^2([0,t_f];\mathbb{R}^m)} .
    % \leq &\ c_7 |t-s| + c_8 |t-s|^{1/2}.
\end{align*}
Since $\{\norm{p_k}_{L^2([0,t_f];\mathbb{R}^m)} \}$ and $\{ \norm{Tp_k}_{C([0,t_f];\mathbb{R}^n)} \}$ both are uniformly bounded for all $p_k \in \mathcal{P}_0$, $\{Tp_k \}$ is uniformly equicontinuous. By the Arzel\`a-Ascoli Theorem \cite{royden2010real}, there is a uniformly convergent subsequence $\{ Tp_{k_j}\} \subset \{Tp_k \}$. 
Without loss of generality, we assume $p_k \rightharpoonup p$ and $Tp_k \rightarrow Tp$ uniformly on $[0,t_f]$, and $\costeval{\eqref{prob: min covariance sensor guidance problem}}{p_k} \leq \costeval{\eqref{prob: min covariance sensor guidance problem}}{p_0}$. We have
\begin{align}
    & \costeval{\eqref{prob: min covariance sensor guidance problem}}{p_0} - \costeval{\eqref{prob: min covariance sensor guidance problem}}{p} \nonumber \\
    = & \costeval{\eqref{prob: min covariance sensor guidance problem}}{p_0} - \costeval{\eqref{prob: min covariance sensor guidance problem}}{p_k} +     \costeval{\eqref{prob: min covariance sensor guidance problem}}{p_k} - \costeval{\eqref{prob: min covariance sensor guidance problem}}{p} \nonumber \\
    \geq & \costeval{\eqref{prob: min covariance sensor guidance problem}}{p_k} - \costeval{\eqref{prob: min covariance sensor guidance problem}}{p}. \label{eq: intermediate step 21}
\end{align}
Hence, to show $\costeval{\eqref{prob: min covariance sensor guidance problem}}{p} \leq \costeval{\eqref{prob: min covariance sensor guidance problem}}{p_0}$, it suffices to show $\costeval{\eqref{prob: min covariance sensor guidance problem}}{p} \leq \liminf_{k \rightarrow \infty} \costeval{\eqref{prob: min covariance sensor guidance problem}}{p_k} $, which is to show
\begin{align}
    & h_f(Tp(t_f)) +\textstyle \int_0^{t_f} h(Tp(t),t) + g(p(t),t) + \trace{\Pi(t)} \dd t \nonumber \\
    \leq & \liminf_{k \rightarrow \infty} h_f(Tp_k(t_f)) + \textstyle  \int_0^{t_f} h(Tp_k(t),t) + g(p_k(t),t) \nonumber \\
    & + \trace{\Pi^k(t)} \dd t ,\label{eq: inequality to be proved}
\end{align}
where $\Pi^k(t)$ is the solution of \eqref{eq: new mild solution of operator Riccati equation} associated with sensor state $Tp_k$. Since $\{Tp_k\}$ converges to $Tp$ uniformly on $[0,t_f]$, the continuity of $h_f(\cdot)$ implies
\begin{equation}
    h_f(Tp(t_f)) = \liminf_{k \rightarrow \infty} h_f(Tp_k(t_f)); \label{eq: intermediate step Fatou's lemma}
\end{equation}
Fatou's lemma \cite{royden2010real} implies
\begin{equation}
     \textstyle \int_0^{t_f} h(Tp(t),t) \dd t \leq \liminf_{k \rightarrow \infty} \textstyle \int_0^{t_f} h(Tp_k(t),t) \dd t;
\end{equation}
and Lemma~\ref{lemma: continuity of trace cost wrt sensor trajectory} implies
\begin{equation}
    \textstyle \int_0^{t_f} \trace{\Pi(t)} \dd t  = \liminf_{k \rightarrow \infty}  \textstyle \int_0^{t_f} \trace{\Pi^k(t)} \dd t. \label{eq: intermediate equality on PDE cost}
\end{equation}
To prove \eqref{eq: inequality to be proved}, based on \eqref{eq: intermediate step Fatou's lemma}--\eqref{eq: intermediate equality on PDE cost}, it suffices to show
\begin{equation}
    \textstyle \int_0^{t_f} g(p(t),t) \dd t \leq \liminf_{k \rightarrow \infty} \textstyle \int_0^{t_f} g(p_k(t),t) \dd t.
\end{equation}
By contradiction, assume there is $\lambda > 0$ such that 
\begin{equation}
    \liminf_{k \rightarrow \infty} \textstyle \int_o^{t_f} g(p_k(t),t) \dd t  < \lambda < \textstyle \int_0^{t_f} g(p(t),t) \dd t. \label{eq: contradiction assumed for existence of solution}
\end{equation}
There exists a subsequence $\{p_{k_j} \} \subset \{p_k \}$ such that ${\{p_{k_j} \} \subset O_{\lambda}  = \{q \in L^2([0,t_f];\mathbb{R}^{m}): \int_0^{t_f} g(q(t),t) \dd t \leq \lambda \}}$.
We wish to show that $O_{\lambda}$ is weakly sequentially closed. By \cite[Theorem~6.1.5]{werner2013optimization}, it suffices to show that $O_{\lambda}$ is convex and closed. Since $g(\cdot,t): \mathbb{R}^{m} \rightarrow \mathbb{R}$ is convex for all $t \in [0,t_f]$, it follows that $O_{\lambda}$ is convex.
Let $\{q_k \} \subset O_{\lambda}$ and $\norm{q_k - q}_{L^2([0,t_f];\mathbb{R}^m)}$ converges to $0$ as $k \rightarrow \infty$. We can choose a subsequence $\{q_{k_j} \} \subset \{q_k  \}$ such that $q_{k_j}$ converges to $q$ pointwise almost everywhere on $[0,t_f]$ \cite[p. 53]{yosida1988functional}. Now we have $g(q_{k_j}(t),t) \geq 0$ for all $t \in [0,t_f]$ (assumption (A11)) and $\lim_{j \rightarrow \infty} g(q_{k_j}(t),t) = g(q(t),t)$ almost everywhere on $[0,t_f]$.
% \begin{enumerate}
%     \item $g(q_{k_j}(t),t) \geq 0$ for all $t \in [0,t_f]$ (assumption (A11));
%     \item $ \lim_{j \rightarrow \infty} g(q_{k_j}(t),t) = g(q(t),t)$ almost everywhere on $[0,t_f]$.
% \end{enumerate}
By Fatou's lemma \cite{royden2010real}, $\int_0^{t_f} g(q(t),t) \dd t  \leq \liminf_{k \rightarrow \infty} \textstyle\int_0^{t_f} g(q_{k_j}(t),t) \dd t \leq \lambda$,
% \begin{equation}
%     \textstyle\int_0^{t_f} g(q(t),t) \dd t  \leq \liminf_{k \rightarrow \infty} \textstyle\int_0^{t_f} g(q_{k_j}(t),t) \dd t \leq \lambda, \label{eq: using Fatou's lemma}
% \end{equation}
where the last inequality holds due to the fact that $\{q_{k_j} \} \subset O_{\lambda}$. Hence, $q \in O_{\lambda}$ and $O_{\lambda}$ is closed.

Since $O_{\lambda}$ is weakly sequentially closed, $p_{k_j} \rightharpoonup p$ implies that $p \in O_{\lambda}$, which contradicts \eqref{eq: contradiction assumed for existence of solution}. Hence, $\costeval{\eqref{prob: min covariance sensor guidance problem}}{p} \leq  \costeval{\eqref{prob: min covariance sensor guidance problem}}{p_0} $ is proved, and we conclude \underline{Condition-2}.

\noindent \textit{\textbf{Proof of \underline{Condition-3}}}: We now show that the mapping $\costeval{\eqref{prob: min covariance sensor guidance problem}}{\cdot}: \mathcal{P} \rightarrow \mathbb{R}$ is weakly sequentially lower semicontinuous on $\mathcal{P}_0$. Suppose $\{p_k \} \subset \mathcal{P}_0$ and $p_k \rightharpoonup p \in \mathcal{P}_0$. We wish to establish $\costeval{\eqref{prob: min covariance sensor guidance problem}}{p} \leq \liminf_{k \rightarrow \infty} \costeval{\eqref{prob: min covariance sensor guidance problem}}{p_k} $, which can be shown using the technique of proving \underline{Condition-2} (starting from \eqref{eq: inequality to be proved}).

Thus we conclude the existence of a solution of problem \eqref{prob: min covariance sensor guidance problem}.
\qed
\end{pf*}

\section{Proof of Theorem~\ref{thm: existence of a solution of the approximated problem}}\label{prf: existence of a solution of the approximated problem}

\vspace{-0.5cm}

\begin{pf*}{Proof}
Since $\trace{\Pi_N(t)} \geq 0$ ($\Pi_N(t)$ is nonnegative and self-adjoint for all $t \in [0,t_f]$) and the mapping $K_N: C([0,t_f];\mathbb{R}^n) \rightarrow \mathbb{R}$ is continuous (see Lemma~\ref{lemma: continuity of trace cost wrt sensor trajectory}), the proof is analogous to that of Theorem~\ref{thm: existence of a solution of minVarGuidanceProblem}, where $\Pi(t)$ is replaced by $\Pi_N(t)$.
\qed
\end{pf*}

\section{Proof of Theorem~\ref{thm: convergence of approximate solution minVarEst}}\label{prf: convergence of approximate solution minVarEst}

\vspace{-0.3cm}

Recall that the notation $\optcosteval{\eqref{prob: approximated min covariance sensor guidance problem}}{p_N^*}$ means the optimal value of \eqref{prob: approximated min covariance sensor guidance problem} evaluated at its optimal solution $p_N^*$, where the dimension of approximation applied to \eqref{prob: approximated min covariance sensor guidance problem} is $N$ (as indicated by the subscript of $p_N^*$). In this section, we attach a subscript to $\eqref{prob: approximated min covariance sensor guidance problem}$, such as $\costeval{\eqref{prob: approximated min covariance sensor guidance problem}_N}{p}$, to indicate its dimension when it is not reflected by the argument, e.g., $\costeval{\eqref{prob: approximated min covariance sensor guidance problem}_N}{p}$ means that the cost of \eqref{prob: approximated min covariance sensor guidance problem} using an $N$-dimensional approximation evaluated at a guidance function $p$. Lemma~\ref{lem: lemma prepared for proving convergence of approximate solution} (see proof in the supplementary material) will be used in the proof of Theorem~\ref{thm: convergence of approximate solution minVarEst}.
\begin{lem}\label{lem: lemma prepared for proving convergence of approximate solution}
    Consider problems \eqref{prob: min covariance sensor guidance problem} and its approximation \eqref{prob: approximated min covariance sensor guidance problem}. Let assumptions (A4)--(A12) hold. Then the following results hold:
    \begin{enumerate}[label={\arabic{enumi}}.,ref={Step \arabic{enumi}},leftmargin=*]
        \item For $p \in C([0,t_f];\mathbb{R}^m)$, $\lim_{N \rightarrow \infty} \costeval{\eqref{prob: approximated min covariance sensor guidance problem}_N}{p} = \costeval{\eqref{prob: min covariance sensor guidance problem}}{p}$;
        \item The mapping $J_{\eqref{prob: min covariance sensor guidance problem}}:C([0,t_f];\mathbb{R}^m) \rightarrow \mathbb{R}^+$ such that $\costeval{\eqref{prob: min covariance sensor guidance problem}}{p} = \int_0^{t_f} \trace{\Pi(t)} \dd t + J_{\text{m}}(\zeta,p)$ is continuous, where $\zeta$ is the sensor state steered by $p$ under the dynamics \eqref{eq: general dynamics of the mobile sensor} and $\Pi(\cdot)$ is the covariance operator obtained through \eqref{eq: new mild solution of operator Riccati equation} with sensor state $\zeta$.
        % \item The mapping $J_{\eqref{prob: approximated min covariance sensor guidance problem}_N}:C([0,t_f];\mathbb{R}^m) \rightarrow \mathbb{R}^+$ such that $\costeval{\eqref{prob: approximated min covariance sensor guidance problem}_N}{p} = \int_0^{t_f} \trace{\Pi_N(t)} \dd t + J_{\text{m}}(\zeta,p)$ is continuous, where $\zeta$ is the sensor state steered by $p$ under the dynamics \eqref{eq: general dynamics of the mobile sensor} and $\Pi_N(\cdot)$ is the approximated covariance operator obtained through \eqref{eq: mild solution of approximate Riccati covariance} with sensor state $\zeta$.
    \end{enumerate}
\end{lem}

\begin{pf*}{Proof of Theorem~\ref{thm: convergence of approximate solution minVarEst}}
        We start with proving \eqref{eq: convergence of the approximate optimal cost}, i.e., $|\optcosteval{\eqref{prob: approximated min covariance sensor guidance problem}}{p_N^*} - \optcosteval{\eqref{prob: min covariance sensor guidance problem}}{p^*}| \rightarrow 0$ as $N \rightarrow \infty$. First,
    \begin{align}
        \optcosteval{\eqref{prob: approximated min covariance sensor guidance problem}}{p_N^*} = &\ \underset{p \in \mathcal{P}(p_{\max},a_{\max})}{\min} \costeval{\eqref{prob: approximated min covariance sensor guidance problem}_N}{p} \nonumber \\
        \leq &\ \costeval{\eqref{prob: approximated min covariance sensor guidance problem}_N}{p^*} \nonumber \\
        \leq &\ |\costeval{\eqref{prob: approximated min covariance sensor guidance problem}_N}{p^*} - \optcosteval{\eqref{prob: min covariance sensor guidance problem}}{p^*}| + \optcosteval{\eqref{prob: min covariance sensor guidance problem}}{p^*}. \nonumber
    \end{align}
    It follows that
    \begin{align}\label{eq: intermediate step 3 for solution convergence}
        \limsup_{N \rightarrow \infty} \optcosteval{\eqref{prob: approximated min covariance sensor guidance problem}}{p_N^*} \leq &\ \optcosteval{\eqref{prob: min covariance sensor guidance problem}}{p^*},
    \end{align}
    because $|\costeval{\eqref{prob: approximated min covariance sensor guidance problem}_N}{p^*} - \optcosteval{\eqref{prob: min covariance sensor guidance problem}}{p^*}| \rightarrow 0$ as $N \rightarrow 0$ (see Lemma~\ref{lem: lemma prepared for proving convergence of approximate solution}-1).
    
    % { \color{red}
    % IMPORTANT: below is the idea of how to prove the result in $\liminf_{N \rightarrow \infty} \optcosteval{\eqref{prob: approximated min covariance sensor guidance problem}}{p_N^*} \geq \optcosteval{\eqref{prob: min covariance sensor guidance problem}}{p^*}$:
    
    % \noindent 1. Choose a convergence subsequence of $\{ \optcosteval{\eqref{prob: approximated min covariance sensor guidance problem}}{p_{N_k}^*} \}_{k=1}^{\infty}$ such that $\optcosteval{\eqref{prob: approximated min covariance sensor guidance problem}}{p_{N_k}^*} \rightarrow \liminf_{N \rightarrow \infty} \optcosteval{\eqref{prob: approximated min covariance sensor guidance problem}}{p_N^*}$ as $k \rightarrow \infty$.
    
    % \noindent 2. Choose a (uniformly) convergent subsequence of $\{p_{N_k} \}_{k=1}^{\infty}$ and denote it by with the same index (to save notation). Denote its limit by $p_{\inf}^*$.
    
    % \noindent 3. Prove $\lim_{k \rightarrow \infty} \optcosteval{\eqref{prob: approximated min covariance sensor guidance problem}}{p_{N_k}^*}  = \costeval{\eqref{prob: min covariance sensor guidance problem}}{p_{\inf}^*} $
    % }
    
    To proceed with proving \eqref{eq: convergence of the approximate optimal cost}, in addition to \eqref{eq: intermediate step 3 for solution convergence}, we shall show $     \liminf_{N \rightarrow \infty} \optcosteval{\eqref{prob: approximated min covariance sensor guidance problem}}{p_N^*} \geq \optcosteval{\eqref{prob: min covariance sensor guidance problem}}{p^*}$.
    Choose a convergent subsequence $\{\optcosteval{\eqref{prob: approximated min covariance sensor guidance problem}}{p_{N_k}^*} \}_{k=1}^{\infty}$ such that 
    \begin{equation}
        \lim_{k \rightarrow \infty} \optcosteval{\eqref{prob: approximated min covariance sensor guidance problem}}{p_{N_k}^*} = \liminf_{N \rightarrow \infty} \optcosteval{\eqref{prob: approximated min covariance sensor guidance problem}}{p_N^*}.
    \end{equation}
    Since the subsequence $\{p_{N_k}^* \}_{k=1}^{\infty} \subset \mathcal{P}(p_{\max},a_{\max})$ which is uniformly equicontinuous and uniformly bounded, by the Arzel\`a-Ascoli Theorem \cite{royden2010real}, there is a (uniformly) convergent subsequence of $\{p_{N_k}^* \}_{k=1}^{\infty}$. We denote this convergent subsequence with the same indices $\{N_k\}_{k=1}^{\infty}$ to simplify notation.
    % Since the guidance functions $p_N^*$ and $p^*$ defined in the set $\mathcal{P}(p_{\max},a_{\max})$ are uniformly equicontinuous and uniformly bounded, by the Arzel\`a-Ascoli Theorem \cite{royden2010real}, there is a uniformly convergent subsequence $\{ p_{N_k}^* \}_{k=1}^{\infty}$ such that $\lim_{k \rightarrow \infty}\optcosteval{\eqref{prob: approximated min covariance sensor guidance problem}}{p_{N_k}^*} = \liminf_{N \rightarrow \infty} \optcosteval{\eqref{prob: approximated min covariance sensor guidance problem}}{p_N^*}$. 
    Denote the limit of $\{p_{N_k}^* \}_{k=1}^{\infty}$ by $p_{\inf}^*$, i.e.,
    \begin{equation}\label{eq: convergence of subsequence of approximate optimal guidance}
        \lim_{k \rightarrow \infty} \norm{p_{N_k}^* - p_{\inf}^*}_{C([0,t_f];\mathbb{R}^m)}  = 0.
    \end{equation}
    
    Next, we show 
    \begin{equation}\label{eq: intermediate step in proving the approximate solution converges}
        \lim_{k \rightarrow \infty} |\optcosteval{\eqref{prob: approximated min covariance sensor guidance problem}}{p_{N_k}^*} - \costeval{\eqref{prob: min covariance sensor guidance problem}}{p_{\inf}^*}| = 0.
    \end{equation}
    First notice that for all $p \in \mathcal{P}(p_{\max},a_{\max})$, $\costeval{\eqref{prob: approximated min covariance sensor guidance problem}_N}{p}$ converges to $\costeval{\eqref{prob: min covariance sensor guidance problem}}{p}$ pointwise as the dimension of approximation $N$ goes to infinity (see Lemma~\ref{lem: lemma prepared for proving convergence of approximate solution}-1). Furthermore, since the sequence of approximated uncertainty cost $\{\int_0^{t_f} \trace{\Pi_N(t)} \dd t \}_{N=1}^{\infty}$ is a monotonically increasing sequence, the sequence $\{ \costeval{\eqref{prob: approximated min covariance sensor guidance problem}_N}{p}\}_{N=1}^{\infty}$ is a monotonically increasing sequence for each $p$ on the compact set $\mathcal{P}(p_{\max},a_{\max})$. By Dini's Theorem \cite[Theorem~7.13]{rudin1964principles}, $|\costeval{\eqref{prob: approximated min covariance sensor guidance problem}_N}{p} - \costeval{\eqref{prob: min covariance sensor guidance problem}}{p}| \rightarrow 0$ uniformly on $\mathcal{P}(p_{\max},a_{\max})$ as $N \rightarrow \infty$. 
    By Moore-Osgood Theorem \cite[Theorem~7.11]{rudin1964principles}, this uniform convergence and the convergence $p_{N_k}^* \rightarrow p_{\inf}^*$ as $k \rightarrow \infty$ (see \eqref{eq: convergence of subsequence of approximate optimal guidance}) imply that
    \begin{align}\label{eq: convergence due to the continuity of original cost function}
        \lim_{k \rightarrow \infty} \costeval{\eqref{prob: min covariance sensor guidance problem}}{p_{N_k}^*} = \lim_{j \rightarrow \infty} \lim_{k \rightarrow \infty} \optcosteval{\eqref{prob: approximated min covariance sensor guidance problem}_j}{p_{N_k}^*}.
    \end{align}
    And the iterated limit in \eqref{eq: convergence due to the continuity of original cost function} equals the double limit \cite[p.~140]{taylor1985general}%https://www.google.com/books/edition/General_Theory_of_Functions_and_Integrat/pczdngEACAAJ?hl=en&gbpv=1&bsq=140
    , i.e.,
    \begin{align}\label{eq: final step to show the convergence of approxiamted subsequence}
        \lim_{j \rightarrow \infty} \lim_{k \rightarrow \infty} \optcosteval{\eqref{prob: approximated min covariance sensor guidance problem}_j}{p_{N_k}^*} = & \lim_{\substack{j \rightarrow \infty \\k \rightarrow \infty}} \optcosteval{\eqref{prob: approximated min covariance sensor guidance problem}_{j}}{p_{N_k}^*} \nonumber \\
        =&  \lim_{k \rightarrow \infty }\optcosteval{\eqref{prob: approximated min covariance sensor guidance problem}}{p_{N_k}^*}.
    \end{align}
    By \eqref{eq: convergence due to the continuity of original cost function}, \eqref{eq: final step to show the convergence of approxiamted subsequence}, and the fact that $\costeval{\eqref{prob: min covariance sensor guidance problem}}{p_{\inf}^*} = \lim_{k \rightarrow \infty} \costeval{\eqref{prob: min covariance sensor guidance problem}}{p_{N_k}^*}$ holds (due to the continuity of $J_{\eqref{prob: min covariance sensor guidance problem}}(\cdot)$, see Lemma~\ref{lem: lemma prepared for proving convergence of approximate solution}-2), we conclude that \eqref{eq: intermediate step in proving the approximate solution converges} holds and
    \begin{align}
        \liminf_{N \rightarrow \infty} \optcosteval{\eqref{prob: approximated min covariance sensor guidance problem}}{p_N^*} = & \  \lim_{k \rightarrow \infty} \optcosteval{\eqref{prob: approximated min covariance sensor guidance problem}}{p_{N_k}^*} \nonumber \\
        = & \ \costeval{\eqref{prob: min covariance sensor guidance problem}}{p_{\inf}^*} \nonumber \\
        \geq & \ \optcosteval{\eqref{prob: min covariance sensor guidance problem}}{p^*}. \label{eq: intermediate step 4 for solution convergence} 
    \end{align}
    Therefore, we conclude $\lim_{N \rightarrow \infty} \optcosteval{\eqref{prob: approximated min covariance sensor guidance problem}}{p_N^*} = \optcosteval{\eqref{prob: min covariance sensor guidance problem}}{p^*}$ from \eqref{eq: intermediate step 3 for solution convergence} and \eqref{eq: intermediate step 4 for solution convergence}.
    % \begin{equation}\label{eq: convergence of the optimal value from approximation to original}
    %     \lim_{N \rightarrow \infty} \optcosteval{\eqref{prob: approximated min covariance sensor guidance problem}}{p_N^*} = \optcosteval{\eqref{prob: min covariance sensor guidance problem}}{p^*}. 
    % \end{equation}
    
    Next, we prove \eqref{eq: convergence of the approximate optimal guidance}, i.e., $|\costeval{\eqref{prob: min covariance sensor guidance problem}}{p_N^*} - \optcosteval{\eqref{prob: min covariance sensor guidance problem}}{p^*}| \rightarrow 0$ as $N \rightarrow \infty$.
    We start with $\optcosteval{\eqref{prob: min covariance sensor guidance problem}}{p^*} \leq \costeval{\eqref{prob: min covariance sensor guidance problem}}{p_N^*}$ for all $N$, which implies that
    \begin{equation}\label{eq: intermediate step 15}
        \optcosteval{\eqref{prob: min covariance sensor guidance problem}}{p^*} \leq \liminf_{N \rightarrow \infty} \costeval{\eqref{prob: min covariance sensor guidance problem}}{p_N^*}.
    \end{equation}
    To prove \eqref{eq: convergence of the approximate optimal guidance}, what remains to be shown is $\optcosteval{\eqref{prob: min covariance sensor guidance problem}}{p^*} \geq \limsup_{N \rightarrow \infty} \costeval{\eqref{prob: min covariance sensor guidance problem}}{p_N^*}$. Choose a convergent subsequence $\{\costeval{\eqref{prob: min covariance sensor guidance problem}}{p_{N_j}^*} \}_{j = 1}^{\infty}$ such that $\lim_{j \rightarrow \infty} \costeval{\eqref{prob: min covariance sensor guidance problem}}{p_{N_j}^*} = \limsup_{N \rightarrow \infty} \costeval{\eqref{prob: min covariance sensor guidance problem}}{p_N^*}$. Since $\{p_{N_j}^* \}_{j=1}^{\infty} \subset \mathcal{P}(p_{\max},a_{\max})$ is uniformly equicontinuous and uniformly bounded, by Arzel\`a-Ascoli Theorem \cite{royden2010real}, $\{p_{N_j}^* \}_{j=1}^{\infty}$ has a (uniformly) convergent subsequence which we denote with the same indices $\{ N_j\}_{j = 1}^{\infty}$ to simplify notation. Denote the limit of $\{p_{N_j}^* \}_{j=1}^{\infty}$ by $p_{\sup}^*$ such that
    \begin{equation}\label{eq: pointwise convergence to p_sup}
        \lim_{j \rightarrow \infty} \norm{p_{N_j}^* - p_{\sup}^*}_{C([0,t_f];\mathbb{R}^m)} = 0.
    \end{equation}
    Due to the continuity of $\costeval{\eqref{prob: min covariance sensor guidance problem}}{\cdot}$ (see Lemma~\ref{lem: lemma prepared for proving convergence of approximate solution}-1), we have $\costeval{\eqref{prob: min covariance sensor guidance problem}}{p_{\sup}^*} = \lim_{j \rightarrow \infty} \costeval{\eqref{prob: min covariance sensor guidance problem}}{p_{N_j}^*} =  \limsup_{N \rightarrow \infty }\costeval{\eqref{prob: min covariance sensor guidance problem}}{p_N^*}$.
    Now we have
    \begin{align}
        & \costeval{\eqref{prob: min covariance sensor guidance problem}}{p_{\sup}^*}  \nonumber \\
        \leq & |\costeval{\eqref{prob: min covariance sensor guidance problem}}{p_{\sup}^*} - \optcosteval{\eqref{prob: min covariance sensor guidance problem}}{p^*}| + \optcosteval{\eqref{prob: min covariance sensor guidance problem}}{p^*} \nonumber \\
        =& |\costeval{\eqref{prob: min covariance sensor guidance problem}}{p_{\sup}^*} - \lim_{N \rightarrow \infty} \optcosteval{\eqref{prob: approximated min covariance sensor guidance problem}}{p_N^*}| + \optcosteval{\eqref{prob: min covariance sensor guidance problem}}{p^*} \nonumber \\
        =& |\costeval{\eqref{prob: min covariance sensor guidance problem}}{p_{\sup}^*} - \lim_{j \rightarrow \infty} \optcosteval{\eqref{prob: approximated min covariance sensor guidance problem}}{p_{N_j}^*}| + \optcosteval{\eqref{prob: min covariance sensor guidance problem}}{p^*}. \label{eq: intermediate step 13}
    \end{align}
    % First notice that for all $p \in \mathcal{P}(p_{\max},a_{\max})$, $\costeval{\eqref{prob: approximated min covariance sensor guidance problem}_N}{p}$ converges to $\costeval{\eqref{prob: min covariance sensor guidance problem}}{p}$ pointwise as the dimension of approximation $N$ goes to infinity (see \eqref{eq: intermediate step 1 for solution convergence}). Furthermore, 
    Since the sequence of approximated uncertainty cost $\{\int_0^{t_f} \Pi_N(t) \dd t \}_{N=1}^{\infty}$ is a monotonically increasing sequence, the sequence $\{ \costeval{\eqref{prob: approximated min covariance sensor guidance problem}_N}{p}\}_{N=1}^{\infty}$ is a monotonically increasing sequence for each $p$ on the compact set $\mathcal{P}(p_{\max},a_{\max})$. Since $\lim_{N \rightarrow \infty}\costeval{\eqref{prob: approximated min covariance sensor guidance problem}_N}{p} = \costeval{\eqref{prob: min covariance sensor guidance problem}}{p}$ for all $p \in \mathcal{P}(p_{\max},a_{\max})$ (see Lemma~\ref{lem: lemma prepared for proving convergence of approximate solution}-1), by Dini's Theorem \cite[Theorem~7.13]{rudin1964principles}, the limit holds uniformly on $\mathcal{P}(p_{\max},a_{\max})$ as $N \rightarrow \infty$. 
    By Moore-Osgood Theorem \cite[Theorem~7.11]{rudin1964principles}, this uniform convergence and the convergence $p_{N_j}^* \rightarrow p_{\sup}^*$ as $j \rightarrow \infty$ (see \eqref{eq: pointwise convergence to p_sup}) imply that $\costeval{\eqref{prob: min covariance sensor guidance problem}}{p_{\sup}^*} = \lim_{k \rightarrow \infty} \lim_{j \rightarrow \infty}  \optcosteval{\eqref{prob: approximated min covariance sensor guidance problem}_k}{p_{N_j}^*}$.
    % \begin{align}\label{eq: convergence due to the continuity of original cost function 2}
    %     \costeval{\eqref{prob: min covariance sensor guidance problem}}{p_{\sup}^*} = \lim_{k \rightarrow \infty} \lim_{j \rightarrow \infty}  \optcosteval{\eqref{prob: approximated min covariance sensor guidance problem}_k}{p_{N_j}^*}.
    % \end{align}
    Furthermore, the iterated limit equals the double limit \cite[p.~140]{taylor1985general}%https://www.google.com/books/edition/General_Theory_of_Functions_and_Integrat/pczdngEACAAJ?hl=en&gbpv=1&bsq=140
    , i.e., 
    \begin{align}
        \lim_{k \rightarrow \infty} \lim_{j \rightarrow \infty}  \optcosteval{\eqref{prob: approximated min covariance sensor guidance problem}_k}{p_{N_j}^*} = & \lim_{\substack{j \rightarrow \infty \\ k \rightarrow \infty}} \optcosteval{\eqref{prob: approximated min covariance sensor guidance problem}_k}{p_{N_j}^*} \nonumber \\
        = & \lim_{j \rightarrow \infty}\optcosteval{\eqref{prob: approximated min covariance sensor guidance problem}}{p_{N_j}^*}.
    \end{align}
    Hence, $\costeval{\eqref{prob: min covariance sensor guidance problem}}{p_{\sup}^*} =  \lim_{j \rightarrow \infty} \optcosteval{\eqref{prob: approximated min covariance sensor guidance problem}}{p_{N_j}^*}$, which, combined with \eqref{eq: intermediate step 13}, implies
    \begin{equation}\label{eq: intermediate step 14}
        \optcosteval{\eqref{prob: min covariance sensor guidance problem}}{p^*} \geq \costeval{\eqref{prob: min covariance sensor guidance problem}}{p_{\sup}^*}  = \limsup_{N \rightarrow \infty} \costeval{\eqref{prob: min covariance sensor guidance problem}}{p_N^*}.
    \end{equation}
    The desired convergence $\lim_{N \rightarrow \infty}\costeval{\eqref{prob: min covariance sensor guidance problem}}{p_N^*} = \optcosteval{\eqref{prob: min covariance sensor guidance problem}}{p^*}$ follows from \eqref{eq: intermediate step 15} and \eqref{eq: intermediate step 14}. \qed
\end{pf*}

% \begin{equation}
% \begin{aligned}
% & \underset{p \in L^2([0,t_f];P)}{\text{minimize}} &&  \int_0^{t_f} \trace{\Pi(t)} \dd t + J_{\text{m}}(\zeta,p)\\
% & \text{subject to}  && \dot{\zeta}(t) = a \zeta(t) + b p(t), \ \zeta(0) = \zeta_0,
% \end{aligned}
% \end{equation}

% \setstretch{1}
\bibliographystyle{abbrv}
\bibliography{reference}

\end{document}

\begin{frontmatter}
\date{01 November 2021}
\title{Supplementary material for: \\
Optimal guidance and estimation of a 2D diffusion-advection process by a team of mobile sensors\thanksref{footnoteinfo}} % Title, preferably not more 

\thanks[footnoteinfo]{This paper was not presented at any IFAC 
meeting. Corresponding author Sheng Cheng.}

\author[Sheng]{Sheng Cheng}\ead{chengs@illinois.edu},    % Add the 
\author[Derek]{Derek A. Paley}\ead{dpaley@umd.edu}               % e-mail address 
% \author[Baiae]{Publius Maro Vergilius}\ead{vergilius@culture.ir}  % (ead) as shown

\address[Sheng]{University of Illinois Urbana-Champaign}  % Please supply                                              
\address[Derek]{University of Maryland, College Park}             % full addresses
% \address[Baiae]{The White House, Baiae}        % here.

\begin{keyword}                           % Five to ten keywords, 
Infinite-dimensional systems; Multi-agent systems; Modeling for control optimization; Guidance navigation and control. 
% chosen from the IFAC 
\end{keyword}                             % keyword list or with the 
                                          % help of the Automatica 
                                          % keyword wizard

\begin{abstract}                          % Abstract of not more than 200 words.
This supplement provides proofs for Lemmas~2.5, 2.7, and C.1 in \cite{Cheng2021estimationInReivew}.
\end{abstract}

\end{frontmatter}

\begin{customlem}{2.5}\label{lemma: continuity of trace cost wrt sensor trajectory}
Let assumptions (A1)--(A3) hold with $q = 1$ and $\Pi \in C([0,t_f];\mathcal{J}_1(\mathcal{H}))$ be defined in \cite[(13)]{Cheng2021estimationInReivew}. If assumption (A4) holds, then the mapping $K(\cdot)$ is continuous.
\end{customlem}

\begin{pf*}{Proof}
Without loss of generality, consider the case of one mobile sensor, i.e., $m_s = 1$. The case of multiple sensors follows naturally. We first show a consequence of the output operator $\mathcal{C}^{\star}$ being continuous with respect to location. Consider two sensor states $\zeta_1, \zeta_2 \in C([0,t_f];\mathbb{R}^n)$. For any $\phi \in \mathcal{H} = L^2(\Omega)$ and all $t \in [0,t_f]$,
\begin{align}
     & \ |\mathcal{C}^{\star}(M\zeta_1(t),t) \phi - \mathcal{C}^{\star} (M\zeta_2(t),t)  \phi| \nonumber \\
     \leq &    \norm{\mathcal{C}(M\zeta_1(t),t) - \mathcal{C}(M \zeta_2(t),t)}_{L^2(\Omega)} \norm{\phi}_{L^2(\Omega)} \nonumber \\
    \leq & \ l\left(|M(\zeta_1(t) - \zeta_2(t))|_2 \right) \norm{\phi}_{L^2(\Omega)}, \label{eq: inequality 1 in proving continuity of PDE cost wrt sensor location} 
\end{align}
where we use the fact that $\mathcal{C}(\cdot,\cdot)$ is the integral kernel of $\mathcal{C}^{\star}(\cdot,\cdot)$. Hence,
\begin{align}\label{eq: intermediate step 1}
    & \norm{\mathcal{C}^{\star}(M\zeta_1(t),t) - \mathcal{C}^{\star} (M\zeta_2(t),t)}_{\mathcal{L}(\mathcal{H};\mathbb{R})} \nonumber \\
    \leq & \ l\left(|M(\zeta_1(t) - \zeta_2(t))|_2 \right).
\end{align}

Since $\mathbb{R}$ is finite-dimensional, there exists $c_1 > 0$ such that\cite[proof of Lemma 4.3]{Burns2015Infinitedimensional}
\begin{align}\label{eq: intermediate step 2}
    &  \norm{\mathcal{C}^{\star}(M\zeta_1(t),t) - \mathcal{C}^{\star} (M\zeta_2(t),t)}_{\mathcal{J}_1(\mathcal{H};\mathbb{R})}\nonumber \\
    \leq &\ c_1 \norm{\mathcal{C}^{\star}(M\zeta_1(t),t) - \mathcal{C}^{\star} (M\zeta_2(t),t)}_{\mathcal{L}(\mathcal{H};\mathbb{R})}.   %\\
    % \leq & c_1 l\left(|M(\zeta_1(t) - \zeta_2(t))|_2 \right).
\end{align}
For brevity, we shall use $\mathcal{C}_1(t)$ for $\mathcal{C}(M\zeta_1(t),t)$ and $\mathcal{C}_2(t)$ for $\mathcal{C}(M\zeta_2(t),t)$. Now,
\begin{align}
    & \norm{\mathcal{C}_1(t) R^{-1} \mathcal{C}_1^{\star}(t) - \mathcal{C}_2(t) R^{-1} \mathcal{C}_2^{\star}(t)}_{\mathcal{J}_1(\mathcal{H})} \nonumber \\
    \leq & \norm{\mathcal{C}_1(t) R^{-1} }_{\mathcal{J}_1(\mathbb{R};\mathcal{H})} \norm{\mathcal{C}_1^{\star} (t) - \mathcal{C}_2^{\star}(t)}_{\mathcal{J}_1(\mathcal{H};\mathbb{R})} \nonumber \\
    & + \norm{R^{-1} \mathcal{C}_2^{\star}(t)}_{\mathcal{J}_1(\mathcal{H};\mathbb{R})} \norm{\mathcal{C}_1(t) - \mathcal{C}_2(t)}_{\mathcal{J}_1(\mathbb{R};\mathcal{H})} \nonumber \\
    =& \left(\norm{\mathcal{C}_1(t) R^{-1} }_{\mathcal{J}_1(\mathbb{R};\mathcal{H})}+ \norm{R^{-1} \mathcal{C}_2^{\star}(t)}_{\mathcal{J}_1(\mathcal{H};\mathbb{R})} \right) \nonumber \\
    & \norm{\mathcal{C}_1^{\star} (t) - \mathcal{C}_2^{\star}(t)}_{\mathcal{J}_1(\mathcal{H};\mathbb{R})}
  \nonumber \\
    \leq & \ c_2  \norm{\mathcal{C}_1^{\star} (t) - \mathcal{C}_2^{\star}(t)}_{\mathcal{J}_1(\mathcal{H};\mathbb{R})} \nonumber \\
    \leq & \ c_2 c_1 l(|M(\zeta_1(t)-\zeta_2(t))|_2) 
    \label{eq: intermediate step 4}
\end{align}
for some $c_2 >0$, where the last inequality follows from \eqref{eq: intermediate step 1} and \eqref{eq: intermediate step 2}. %(The equality in the third line holds roughly because for $A \in \mathcal{L}(X,Y) $, we have $\text{dom}(A^{\star}) = Y^{\star}$ and $\norm{A^{\star}}_{\mathcal{L}(Y^{\star},X^{\star})} = \norm{A}_{\mathcal{L}(X,Y)}$).

By \cite[(13)]{Cheng2021estimationInReivew}, the mapping $\Pi: [0,t_f] \rightarrow \mathcal{J}_1(\mathcal{H})$ varies continuously in $\sup_{t \in [0,t_f]} \norm{\cdot}_{\mathcal{J}_1(\mathcal{H})}$-norm with respect to $\bar{\mathcal{C}} \bar{\mathcal{C}}^{\star}(\cdot)$~\cite{Burns2015Infinitedimensional}. Hence, there exists $c_3 >0$ such that 
\begin{align}
     & \sup_{t\in [0,t_f]} \norm{\Pi_{1}(t) - \Pi_{2}(t)}_{\mathcal{J}_1(\mathcal{H})}  \nonumber \\
    \leq & \sup_{t\in [0,t_f]} c_3 \norm{\bar{\mathcal{C}}_1 \bar{\mathcal{C}}_1^{\star}(t) - \bar{\mathcal{C}}_2 \bar{\mathcal{C}}_2^{\star}(t)}_{\mathcal{J}_1(\mathcal{H})}. \label{eq: intermediate step 12}
\end{align}

Now, we have
\begin{align}
& |K(\zeta_1)- K(\zeta_2)| \nonumber \\
    % = & \left| \int_0^{t_f} \trace{\Pi_1(t)} - \trace{\Pi_2(t)} \dd t \right| \nonumber \\
    % = & \left| \int_0^{t_f} \norm{\Pi_1(t)}_{\mathcal{J}_1(\mathcal{H})} - \norm{\Pi_2(t)}_{\mathcal{J}_1(\mathcal{H})} \dd t \right| \nonumber \\
    = & \left|\int_0^{t_f}   \norm{\Pi_1(t)}_{\mathcal{J}_1(\mathcal{H})} - \norm{\Pi_2(t)}_{\mathcal{J}_1(\mathcal{H})} \dd t \right|\nonumber \\
    \leq & \int_0^{t_f} \left|  \norm{\Pi_1(t)}_{\mathcal{J}_1(\mathcal{H})} - \norm{\Pi_2(t)}_{\mathcal{J}_1(\mathcal{H})} \right|\dd t \nonumber \\
    \leq & \int_0^{t_f}   \norm{\Pi_1(t) - \Pi_2(t)}_{\mathcal{J}_1(\mathcal{H})}\dd t \nonumber \\
    \leq & \sup_{t \in [0,t_f]}\norm{\Pi_1(t) - \Pi_2(t)}_{\mathcal{J}_1(\mathcal{H})}t_f. \label{eq: intermediate step 16}
\end{align}
It follows from \eqref{eq: intermediate step 4}--\eqref{eq: intermediate step 16} that
\begin{equation*}
    |K(\zeta_1)- K(\zeta_2)| \leq c_1 c_2 c_3 t_f  \sup_{t\in [0,t_f]}   l(|M(\zeta_1(t)-\zeta_2(t))|_2) ,
\end{equation*}
and we conclude the continuity of $K(\cdot)$. \qed
\end{pf*}

\begin{customlem}{2.7}\label{lem: continuity of approximated trace cost wrt sensor trajecotry}
Let assumptions (A5)--(A7) hold and $\Pi_{N}(t)$ be defined as in \cite[(16)]{Cheng2021estimationInReivew}. If assumption (A8) holds, then the mapping $K_N: C([0,t_f];\mathbb{R}^n) \rightarrow \mathbb{R}^+$ such that $K_N(\zeta) \defeq \int_0^{t_f} \trace{\Pi_N(t)} \dd t$ is continuous.
\end{customlem}

\begin{pf*}{Proof}
    Since the norm defined on $\mathcal{H}_N$ is inherited from that of $\mathcal{H}$, the proof follows from the derivation of Lemma~\ref{lemma: continuity of trace cost wrt sensor trajectory}'s proof. \qed
\end{pf*}

\begin{customlem}{C.1}\label{lem: lemma prepared for proving convergence of approximate solution}
    Consider problems \cite[(P)]{Cheng2021estimationInReivew} and its approximation \cite[(AP)]{Cheng2021estimationInReivew}. Let assumptions (A4)--(A12) hold. Then the following results hold:
    \begin{enumerate}[label={\arabic{enumi}}.,ref={Step \arabic{enumi}},leftmargin=*]
        \item For $p \in C([0,t_f];\mathbb{R}^m)$, $\lim_{N \rightarrow \infty} \costeval{\rm{(AP)}_N}{p} = \costeval{\rm{(P)}}{p}$;
        \item The mapping $J_{\rm{(P)}}:C([0,t_f];\mathbb{R}^m) \rightarrow \mathbb{R}^+$ such that $\costeval{\rm{(P)}}{p} = \int_0^{t_f} \trace{\Pi(t)} \dd t + J_{\text{m}}(\zeta,p)$ is continuous, where $\zeta$ is the sensor state steered by $p$ under the dynamics \cite[(2)]{Cheng2021estimationInReivew} and $\Pi(\cdot)$ is the covariance operator obtained through \cite[(13)]{Cheng2021estimationInReivew} with sensor state $\zeta$.
        % \item The mapping $J_{\rm{(AP)}_N}:C([0,t_f];\mathbb{R}^m) \rightarrow \mathbb{R}^+$ such that $\costeval{\rm{(AP)}_N}{p} = \int_0^{t_f} \trace{\Pi_N(t)} \dd t + J_{\text{m}}(\zeta,p)$ is continuous, where $\zeta$ is the sensor state steered by $p$ under the dynamics \eqref{eq: general dynamics of the mobile sensor} and $\Pi_N(\cdot)$ is the approximated covariance operator obtained through \eqref{eq: mild solution of approximate Riccati covariance} with sensor state $\zeta$.
    \end{enumerate}
\end{customlem}

\begin{pf*}{Proof of Lemma~\ref{lem: lemma prepared for proving convergence of approximate solution}}
    
    1. We first prove that for $p \in C([0,t_f];\mathbb{R}^m)$, %$p \in \mathcal{P}(p_{\max},a_{\max})$, 
    \begin{equation}\label{eq: intermediate step 1 for solution convergence}
    \lim_{N \rightarrow \infty} |\costeval{\rm{(AP)}_N}{p} - \costeval{\rm{(P)}}{p}| = 0.
    \end{equation}
    
    To establish \eqref{eq: intermediate step 1 for solution convergence}, it suffices to show
    \begin{equation}
        \lim_{N \rightarrow \infty} \left| \int_0^{t_f} \trace{\Pi_N(t)} - \trace{\Pi(t)} \dd t \right| = 0.
    \end{equation}
    We have
    \begin{align}
    & \left| \int_0^{t_f} \trace{\Pi_N(t)} - \trace{\Pi(t)} \dd t \right| \nonumber \\
    = & \left| \int_0^{t_f} \lVert \Pi_N(t) \rVert_{\mathcal{J}_1(\mathcal{H})} - \lVert \Pi(t) \rVert_{\mathcal{J}_1(\mathcal{H})} \text{d} t \right| \nonumber \\
    \leq & \int_0^{t_f} \left|  \norm{\Pi_N(t)}_{\mathcal{J}_1(\mathcal{H})} - \norm{\Pi(t)}_{\mathcal{J}_1(\mathcal{H})} \right|\dd t \nonumber \\
    \leq & \int_0^{t_f}   \norm{\Pi_N(t) - \Pi(t)}_{\mathcal{J}_1(\mathcal{H})}\dd t \nonumber \\
    \leq & \sup_{t \in [0,t_f]}\norm{\Pi_N(t) - \Pi(t)}_{\mathcal{J}_1(\mathcal{H})}t_f. \label{eq: arguing convergence of trace of covariance operator}
\end{align}
     By \cite[Theorem~2.6]{Cheng2021estimationInReivew}, $\sup_{t \in [0,t_f]} \norm{\Pi_N(t) - \Pi(t)}_{\mathcal{J}_q(\mathcal{H})} \rightarrow 0$ as $N \rightarrow \infty$. And specifically, when $q = 1$, the convergence in \eqref{eq: intermediate step 1 for solution convergence} holds due to \eqref{eq: arguing convergence of trace of covariance operator}.
     
     2. The cost function of $\rm{(P)}$ is the sum of two parts: the uncertainty cost $\int_0^{t_f} \trace{\Pi(t)} \dd t$, cast as a continuous mapping $K: C([0,t_f];\mathbb{R}^n) \rightarrow \mathbb{R}^+$ (see Lemma~\ref{lemma: continuity of trace cost wrt sensor trajectory}) and the mobility cost $J_{\text{m}}(\zeta,p)$, cast as a mapping $\bar{J}_{\text{m}}: C([0,t_f];\mathbb{R}^m) \rightarrow \mathbb{R}^+$  which we define below. The mapping $\bar{J}_{\text{m}}$ is the single argument version of the original mobility cost by defining the sensor state as a mapping of the sensor guidance. \revision{Here, we redefine the domain of the map $T$ in the proof of \cite[Theorem~3.1]{Cheng2021estimationInReivew} such that $T:C([0,t_f];\mathbb{R}^m) \rightarrow C([0,t_f];\mathbb{R}^n)$. The continuity of $T$ still holds \cite{werner2013optimization}, i.e., for $p_1,p_2 \in C([0,t_f];\mathbb{R}^m)$ there exist $c_5 >0$ such that
    \begin{align}
        \norm{Tp_1 - Tp_2}_{C([0,t_f];\mathbb{R}^n)} \leq  c_5 \norm{p_1 - p_2}_{C([0,t_f];\mathbb{R}^m)}. \label{eq: boundedness of trajectory norm}
    \end{align}} 
% The mapping $T$ is continuous because for $p_1,p_2 \in C([0,t_f];\mathbb{R}^m)$
%     \begin{align}
%         & |(Tp_1)(t) - (Tp_2)(t)|_1 \nonumber \\
%         %= & \int_0^t e^{\alpha (t-\tau)} \beta (p_1(\tau)-p_2(\tau)) \dd \tau \nonumber \\
%         \leq & \int_0^t |e^{\alpha (t - \tau)} \beta |_1 \dd \tau \sup_{\tau \in [0,t]} |p_1(\tau) - p_2(\tau)|_1,
%     \end{align}
%     which implies the existence of $c_5 >0$ such that
%     \begin{align}
%         \norm{Tp_1 - Tp_2}_{C([0,t_f];\mathbb{R}^n)} \leq  c_5 \norm{p_1 - p_2}_{C([0,t_f];\mathbb{R}^m)}. \label{eq: boundedness of trajectory norm}
%     \end{align}
    
    \noindent Let $\bar{J}_{\text{m}}(p) \defeq J_{\text{m}}(Tp,p)$ and we show $\bar{J}_{\text{m}}$ is continuous.
    Define mappings $G: C([0,t_f];\mathbb{R}^m) \rightarrow \mathbb{R}^+$, $H: C([0,t_f];\mathbb{R}^n) \rightarrow \mathbb{R}^+$, and $H_f: C([0,t_f];\mathbb{R}^n) \rightarrow \mathbb{R}^+$ such that
    \begin{align}
        G(p) = & \ \int_0^{t_f} g(p(t),t) \dd t, \\
        H(p) = & \ \int_0^{t_f} h(Tp(t),t) \dd t, \\
        H_f(p) = & \  h_f(Tp(t_f)).
    \end{align}
    Since $\bar{J}_{\text{m}}(p) = G(p) + H(p) + H_f(p)$, we shall proceed with showing that the mappings $G$, $H$, and $H_f$ are continuous.

    Let $p_1,p_2 \in \mathcal{P}(p_{\max},a_{\max})$. Both the set of admissible guidance's values $P_0 \defeq \cup_{t \in [0,t_f]} \{p(t): p \in \mathcal{P}(p_{\max},a_{\max}) \}$ and the interval $[0,t_f]$ are closed and bounded (hence compact). Since $g: P_0 \times [0,t_f] \rightarrow \mathbb{R}^+$ is continuous, by the Heine-Cantor Theorem \cite[Proposition~5.8.2]{sutherland2009introduction}, $g$ is uniformly continuous, i.e., for all $\epsilon > 0$ there exists $\delta > 0$ such that for all $t \in [0,t_f]$, $|p_1(t)-p_2(t)| < \delta$ implies $|g(p_1(t),t) - g(p_2(t),t)| < \epsilon$.
    Hence, it follows that
    \begin{gather}
        \norm{p_1-p_2}_{C([0,t_f];\mathbb{R}^m)} =  \sup_{t \in [0,t_f]} |p_1(t)-p_2(t)| < \delta, \nonumber \\
        \Rightarrow |g(p_1(t),t) - g(p_2(t),t)| <  \epsilon, \quad \forall t \in [0,t_f].
    \end{gather}
    Therefore, for all $\epsilon>0$ there exists $\delta >0$ such that $\norm{p_1-p_2}_{C([0,t_f];\mathbb{R}^m)} < \delta $ implies
    \begin{equation}
        \int_0^{t_f} |g(p_1(t),t) - g(p_2(t),t)|  \dd t < \epsilon t_f,
    \end{equation}
    which establishes the continuity of the mapping $G$.

    Since the continuous image of a compact set is compact \cite[Proposition~5.5.1]{sutherland2009introduction}, the image set $T(\mathcal{P}(p_{\max},a_{\max}))$ is compact, i.e., the set $\Xi \defeq \{\zeta \in C([0,t_f];\mathbb{R}^n): \zeta = Tp, p \in \mathcal{P}(p_{\max},a_{\max}) \}$ is compact. The compactness of $\Xi$ implies that the set of sensor state's values $\zeta(t)$, $\Xi_0 \defeq \cup_{t \in [0,t_f]} \{\zeta(t)|\zeta \in \Xi \}$, is closed. Furthermore, since $\norm{Tp}_{C([0,t_f];\mathbb{R}^n)}$ is bounded (see \eqref{eq: boundedness of trajectory norm}) and $\Xi_0$ is finite dimensional, the set $\Xi_0$ is compact. The compactness of $\Xi_0$ and continuity of the function $h: \Xi_0 \times [0,t_f] \rightarrow \mathbb{R}^+$ implies that $ h$ is uniformly continuous by the Heine-Cantor Theorem \cite[Proposition~5.8.2]{sutherland2009introduction}. Hence, for all $\epsilon > 0$ there exists $\delta > 0$ such that if ${\norm{p_1-p_2}_{C([0,t_f];\mathbb{R}^m)} < \delta/c_5}$, which implies $\norm{Tp_1-Tp_2}_{C([0,t_f];\mathbb{R}^n)} < \delta$, then
    \begin{equation}
         \int_0^{t_f} |h(Tp_1(t),t) - h(Tp_2(t),t)| \dd t \leq \epsilon t_f,
    \end{equation}
    which concludes the continuity of the mapping $H$.
    
    The mapping $H_f$ is continuous because for all $\epsilon > 0$ there exists $\delta >0$ such that if $\norm{p_1 - p_2}_{C([0,t_f];\mathbb{R}^m)} < \delta/c_5 $, which implies $\sup_{t \in [0,t_f]} |Tp_1(t) - Tp_2(t)|< \delta$, then 
    \begin{equation}
        |Tp_1(t_f) - Tp_2(t_f)| < \delta.
    \end{equation}
    Furthermore, $|H_f(p_1) - H_f(p_2)| = |h_f(Tp_1(t_f)) - h_f(Tp_2(t_f))| < \epsilon$ holds due to the continuity of $h_f$.
    
    Hence, we conclude the continuity of $\bar{J}_{\text{m}}$, which, together with the continuity of $K(\cdot)$ and \eqref{eq: boundedness of trajectory norm}, implies the continuity of $J_{\rm{(P)}}(\cdot)$. \qed
    % 3. The proof of the continuity of $J_{\rm{(AP)}}(\cdot)$ is identical to that of $J_{$\rm{(P)}$}(\cdot)$, except that $\Pi(\cdot)$ and $K(\cdot)$ are replaced by $\Pi_N(\cdot)$ and $K_N(\cdot)$ (see Lemma~\ref{lem: continuity of approximated trace cost wrt sensor trajecotry}), respectively.
\end{pf*}

\bibliographystyle{abbrv}
\bibliography{reference}